\documentclass{article}
\title{Comparison of Numerical Solvers for Differential Equations for Holonomic Gradient Method in Statistics}
\author{N.Takayama, T.Yaguchi, Y.Zhang}
\date{April 9, 2026}
\usepackage{a4wide}
\usepackage{color}
\usepackage{url}
\usepackage{hyperref}      
\usepackage{graphicx}

\newtheorem{example}{Example}
\newtheorem{algorithm}{Algorithm}
\newtheorem{definition}{Definition}
\newtheorem{lemma}{Lemma}
\newtheorem{proposition}{Proposition}
\newtheorem{theorem}{Theorem}
\newtheorem{situation}{Situation}

\def\onProg{Programs of this paper can be obtainable
from \\ {\color{blue} \url{https://www.math.kobe-u.ac.jp/OpenXM/Math/defusing/ref.html}}.
The folder/program names are given in blue letters.
(Ignore the folder names if the folder does not exist.)
{\color{magenta} The letter ``-'' in the file name may be ``\_'' for some cases.}
}
\def\prog#1{({\color{blue} #1})}
\def\progweb#1{({\color{blue} #1})}
\def\prog#1{ }
\def\progweb#1{ }
\def\onProg#1{}

\def\qed{\hfill $//$ \bigbreak}
\def\F0true{F_0^{\rm true}}
\def\pd#1{\partial_{#1}}
\def\RR{\textbf{R}}

\begin{document}
\maketitle

\section{Introduction}
Definite integrals with parameters of holonomic functions satisfy
holonomic systems of linear partial differential equations.
When we restrict parameters to a one dimensional curve,
then the system becomes a linear ordinary differential equation (ODE).
For example, the definite integral with a parameter $t$
$$
  Z(t) = \int_0^{+\infty} \exp(ty-y^3)dy
$$
satisfies the ODE $ (3 \pd{t}^2-t) Z(t) = 1 $
where $\pd{t}=\frac{d}{dt}$ (see \cite[Chap 6]{dojo-en} 
as an introductory exposition).
We can evaluate the integral by solving the linear ODE numerically.
The method of solving the linear ODE
is not a high precision method, but it gives values of the integral
at several evaluation points at once.
This approach to evaluate numerically definite integrals is called
{\it the holonomic gradient method} (HGM)
and it is useful to evaluate several normalizing constants 
in statistics. See references in \cite{hgm} 
or in ,e.g., \cite{HNTT}, \cite{TJKZ}.

The HGM consists of 3 steps.
The first step is to find a linear ODE for the definite integral $Z(t)$ with
a parameter $t$. 
The 2nd step is to find value of $Z(t)$ and its derivatives at some points
of $t$.
The 3rd step is to solve the ODE.
We can apply several results in computer algebra and in the theory of
special functions of several variables for the first and 2nd steps.
Results in numerical integration can be applied to the 2nd step.
The HGM is successful for a broad class of distributions and
the 3rd step is relatively easy in most of them.
When the function $Z(t)$ is the dominant solution
(the maximal growth rate solution)
for the ODE when $t \rightarrow +\infty$ as in the case
\cite[Theorem 2]{HNTT}, the 3rd step can be performed by a simple application
of standard numerical solver for ODEs.
However, there are some HGM problems for which we encounter difficulties
in the 3rd step.
The difficulty is that
when the function $Z(t)$ is not dominant and 
we solve the linear ODE in a long interval,
the numerical solution might become a false one.
We have this difficulty
in our analysis of the outage probability
of MIMO WiFi systems (see the section \ref{sec:Hkn} and \cite{DOTS})
and in solving an ODE numerically of the expectation of the Euler characteristic
of a random manifold (see the section \ref{sec:euler} and \cite{TJKZ}).
In principle, if we give an initial value vector exactly and 
the time step is sufficiently small, we do not get a false solution.
However, these assumptions are not always satisfied in practice.

In this paper, we will discuss and compare 
methods to solve linear ODE's in the 3rd step of the HGM.
We intend to provide references to the HGM users
for choosing numerical solvers of ODE's for the HGM.
We will present methods that work and methods that does not work.
Of course, it depends on a given problem.

\section{Methods to test}  \label{sec:methods}
We will compare some methods mainly from a point of view
of using for the 3rd step of the HGM\footnote{\onProg}.
In the HGM, the function $Z(t)$ to evaluate has an integral representation
and satisfies an ODE 
\begin{equation} \label{eq:ode}
 Lf=b, \quad L=\sum_{k=0}^r c_k(t) \pd{t}^k, \ \pd{t}=\frac{d}{dt}
\end{equation}
where $b(t)$ is an inhomogeneous term and $f(t)$ is an unknown function.
See, e.g.,  \cite{DOTS}, \cite{HNTT}, and papers in \cite{hgm} as to examples.
The differential equation (\ref{eq:ode}) might have a solution $u$
which dominates $Z$ at $+\infty$. 
In other words, there might be a dominant solution $u$ such that
$|u/Z| \rightarrow \infty$ when $t \rightarrow +\infty$.
In such cases, it is not easy to obtain numerical values of $Z(t)$ globally
by solving the ODE even when the almost accurate initial values of 
$Z(t), Z'(t), \ldots$ are known. 
See Examples \ref{ex:easy} and \ref{ex:airy1}.
We call such case {\it the unstable HGM problem}\/.
We want to obtain correct values of $Z$ as large interval of $t$ as possible.

Here is a list of methods we try.
\begin{enumerate}
\item The Runge-Kutta method (see, e.g., \cite[Chapter 2]{HNW}). \label{item:RK}
\item The implicit Runge-Kutta method (see, e.g., \cite[Chapters 2 and 4]{HNW}). \label{item:iRK}
\item Multi-step methods (see, e.g., \cite[Chapters 3 and 5]{HNW}). \label{item:BDF}
\item The discrete QR method, the continuous orthonormalization, the Riccati transformation \cite{conte-1966}, \cite{AMR}, \cite{Davey1983}, \cite{Acton}, \cite{DOR1988}, \cite{DRV-1994}. 
\item The spectral method in the approximation theory
(see, e.g., \cite{ApproxFun}, \cite{OT}, \cite{Trefethen}, Section \ref{sec:chebyshef}). \label{item:SM}
\item Defusing method (filter method): restricting the initial value vector to a sub linear space
(see, e.g., Section \ref{sec:defusing}). \label{item:defusing}
\item Sparse interpolations/extrapolations by ODE: solving a generalized boundary value problem 
by solving a linear equation or by solving an optimization
problem (see, e.g., \cite{AMR}, \cite{boyd-2001}, \cite{canuto-2006} and the Section \ref{sec:sparse-interpolation}). \label{item:boundary}
\end{enumerate}
We compare these methods by different implementations on different languages:
\cite{ApproxFun}, \cite{chebfun}, \cite{deSolve}, \cite{gsl}, \cite{MathematicaRK},
\cite{mpfr}, \cite{scipy_linalg}, \cite{solve_ivp}, \cite{least_square},
\cite{cvxopt_qp}. 
The robustness for input errors is an important point in our comparison,
because initial values or boundary values might be evaluated
by Monte-Carlo integrators.

The first 5 approaches are standard.
In the HGM, we may be able to evaluate the target function $Z(t)$
at any point by numerical integration methods. 
Of course, the value obtained is not always accurate.
For example, the accuracy of Monte Carlo integrators are not high.
In this case, Step 3 of the HGM can be completed 
by solving the interpolation problem or the generalized boundary value problem.
On the other hand, while approximating the integral is difficult, 
the value at a certain point can sometimes be calculated in detail 
using a series (see, e.g., \cite{HNTT}.
In this case, the initial value problem should be solved in Step 3 of the HGM.

The last two methods of the defusing method and sparse interpolation/extrapolation method by ODE are expected to be used mainly under these situations of the HGM.

The program codes of this paper are obtainable from
\url{https://www.math.kobe-u.ac.jp/OpenXM/Math/defusing/ref.html}
\cite{web-defusing}.

\section{Defusing method (filter method) and discrete QR method} \label{sec:defusing}

In the HGM, we do not always have an exact initial value vector. 
We often have the case that the target function (the target integral) has a moderate growth, but the ODE has rapid growth solutions, too. Then, we need to correct the initial value vector so that the rapid growth solution does not cause a blow-up of the numerical solution by an error of the initial value vector. 

This instability is a classic hallmark of stiff differential equations. 
In the context of boundary-value problems, the problem of dominant solutions overwhelming subdominant ones has historically been addressed by techniques such as continuous orthonormalization \cite{Davey1983},
discrete QR decomposition method \cite{DRV-1994}, 
or the Riccati transformation method \cite{DOR1988}. 
Continuous orthonormalization, for instance, integrates forward while periodically applying Gram-Schmidt orthogonalization to prevent the solution vectors from becoming linearly dependent on the dominant growing mode. 
This method works well even when there are numerical errors in boundary values. \progweb{orthogonal/davey17-bv.rr}
This method can also be used for initial value problems, but adaptive step control is difficult when explosive solutions are involved.
\progweb{orthogonal/davey11-airy.rr}

The discrete QR decomposition method is also a powerful method for stably finding solutions of boundary value problems (see, e.g., \cite{conte-1966}, \cite{AMR}). 
While the fundamental solution matrix Q, obtained as an orthogonal matrix, can be calculated stably, the triangular matrix R used to correct it and obtain the solution to the original equation honestly amplifies the error in the initial values, making it difficult to find a subdominant solution as the solution to the initial value problem. \progweb{orthogonal/davey16-pkg.rr}

While these are popular stabilization techniques for boundary value problems, 
we propose a discrete, heuristic projection approach, which we will call the \textit{defusing method}, to avoid a blow-up of a solution 
of an initial value problem under the specific situation of the HGM. 
Any solution of the ODE can be expressed as a linear combination of basis functions of different growth orders. 
The defusing method removes some basis functions of specified growth order to obtain the correct solution. 
This method is an analogy of filtering some specified frequencies from wave data, effectively projecting the erroneous initial value vector onto the stable subspace spanned by the subdominant eigenvectors.
This method can be used for solving initial value problems 
with conventional Runge-Kutta shemes and also with
the discrete QR decomposition method.
This is likely a well-known empirical method among experts, but since it can be conveniently used with HGM, we will explain this method below.

In the HGM, we want to evaluate numerically 
a definite integral with a parameter $t$ on some interval of $t$.
We put this function $f(t)$.
The HGM gives an algebraic method to find an ordinary differential equation
(ODE)
for $F(t)=(f(t),f'(t), f''(t), \ldots, f^{(r-1)})^T$
wher $r$ is the rank of the ODE.
We evaluate the integral by Monte-Carlo integration methods
or by the saddle point approximation for large $t$'s.
Since $F(t)$ satisfies the differential equation,
we can apply algorithms 
(see, e.g., \cite{barkatou1997}, \cite{hoeij1997} )
to find a basis of series solutions at $t=\infty$.
Let $F_{*j}$, $j=1, \ldots, r$ be the basis of series solutions around
$t=\infty$ of the ODE.
We denote by $F_{ij}$ is the $i$th component of $F_{*j}$
and $F_i$ is the $i$th component of $F$.
Any solution of the ODE is written as a linear combination of $F_{*j}$'s
over ${\bf C}$.
We assume 
$$ F_{11} > F_{12} > \cdots > F_{1m} > \cdots > F_{1r} > 0
$$
with different growth orders
when $t \rightarrow +\infty$.
We suppose that
\begin{equation} \label{eq:find_m}
  \left| \frac{\mbox{numerical approximation of $F_1(t)$}}
                {F_{1j}(t)}
    \right|
  \rightarrow 0, \quad  t \rightarrow +\infty
\end{equation}
for $ 1 \leq j < m$.
See (\ref{eq:find_m_Hkn}) as an example.
Under these assumptions, we can conclude that
$F(t)$ is approximated as a linear combination of
$F_{*m}, \ldots, F_{*r}$.
The coefficients of the linear combination can be estimated 
by numerical approximation values of $f(t)$ for large $t$'s.
Thus, we can interpolate and extrapolate numerical values  of $F(t)$
around $t=+\infty$
by a set of fundamental solutions of the ODE.
We call this method the defusing method by series solutions
around $t=+\infty$. 
See the Example \ref{ex:hkn_asymp}.

We want to find a numerical solution of the initial value problem 
of the ordinary differential equation (ODE)
\begin{eqnarray}
 \frac{dF}{dt} &=& P(t) F   \label{eq:ode1} \\
 F(t_0) &=& \F0true  \in {\bf R}^r \label{eq:init1}
\end{eqnarray}
where $P(t)$ is an $r \times r$ matrix,
$F(t)$ is a column vector function of size $r$,
and $\F0true $ is the initial value of $F$ at $t=t_0$.
Let us explain a defusing method for the initial value problem.
Solving this problem is the final step of the holonomic gradient method (HGM)
\cite{hgm}.
We often encounter the following situation in the final step.
\begin{situation}\rm \label{situation123}
\begin{enumerate}
\item The initial value has at most 3 digits of accuracy. We denote this initial value
$F_0$.
\item The property $|F| \rightarrow 0 $ when $t \rightarrow +\infty$ is known, e.g., 
from a background of the statistics.
\item There exists a solution ${\tilde F}$ of (\ref{eq:ode1}) such that
$|{\tilde F}| \rightarrow +\infty$ or non-zero finite value
when $t \rightarrow +\infty$.
\end{enumerate}
\end{situation}
Under this situation, the HGM works only for a very short interval of $t$
because the error of the initial value vector makes the fake solution 
${\tilde F}$ dominant and it hides the true solution $F(t)$.
We call this bad behavior of the HGM {\it the instability of the HGM}\/.

\begin{example}  \rm \label{ex:easy}
$$
\frac{d}{dt} F = \left(
\begin{array}{ccc}
-1 & 1 & 0 \\
 0 & -1 & 1 \\
 0 & 0 & 0
\end{array}
\right) F
$$
The solution space is spanned by
$F^1=(\exp(-t),0,0)^T$, $F^2 = (0,\exp(-t),0)^T$,
$F^3=(1,1,1)^T$.
The initial value $(1,0,0)^T$ at $t=0$ yields the solution $F_1$.
Add some errors 
$(1,10^{-30},10^{-30})^T$
to the initial value.
Then, we have
\begin{center}
\begin{tabular}{|l|l|l|}
$t$ & value $F_1$ by RK & difference $F_1-F^1_1$ \\ \hline
50 & 1.92827e-22 & 9.99959e-31 \\
60 & 8.75556e-27 & 1.00000e-30 \\
70 & 1.39737e-30 & \fbox{1.00000e-30} \\
80 & 1.00002e-30 & \fbox{1.00000e-30} \\
\end{tabular}
\end{center}
We can see the instability.
\end{example}

\begin{example} \rm \label{ex:airy1}
$$ P(t) = \left(\begin{array}{cc} 0 & 1 \\ t & 0 \\ \end{array}\right). $$
This differential equation is obtained from the Airy differential equation
$$ y''(t) - t y(t) = 0 $$
by putting $F=(y(t),y'(t))^T$.
It is well-known that the Airy function
$$ {\rm Ai}(t) = \frac{1}{\pi} \lim_{b \rightarrow +\infty}
\int_0^b  \cos\left(\frac{s^3}{3}+ts \right) ds
$$
is a solution of the Airy differential equation and
\begin{eqnarray*}
 {\rm Ai}(0) &=& \frac{1}{3^{2/3} \Gamma(2/3)} = 0.355028053887817\cdots \\
 {\rm Ai}'(0) &=& \frac{1}{3^{1/3} \Gamma(1/3)} = -0.258819403792807\cdots \\
 \lim_{t \rightarrow +\infty} {\rm Ai}(t) &=& 0 \\
 \lim_{t \rightarrow +\infty} {\rm Ai}'(t) &=& 0 \\
\end{eqnarray*}
Figure \ref{fig:airy} is a graph of Airy Ai function and Airy Bi function.
The function $F(t) = ({\rm Ai}(t), {\rm Ai}'(t))^T$ satisfies 
the condition 2 of the Situation \ref{situation123} of the instability problem.

\begin{figure}[htb]
\begin{center}
\includegraphics[width=7cm]{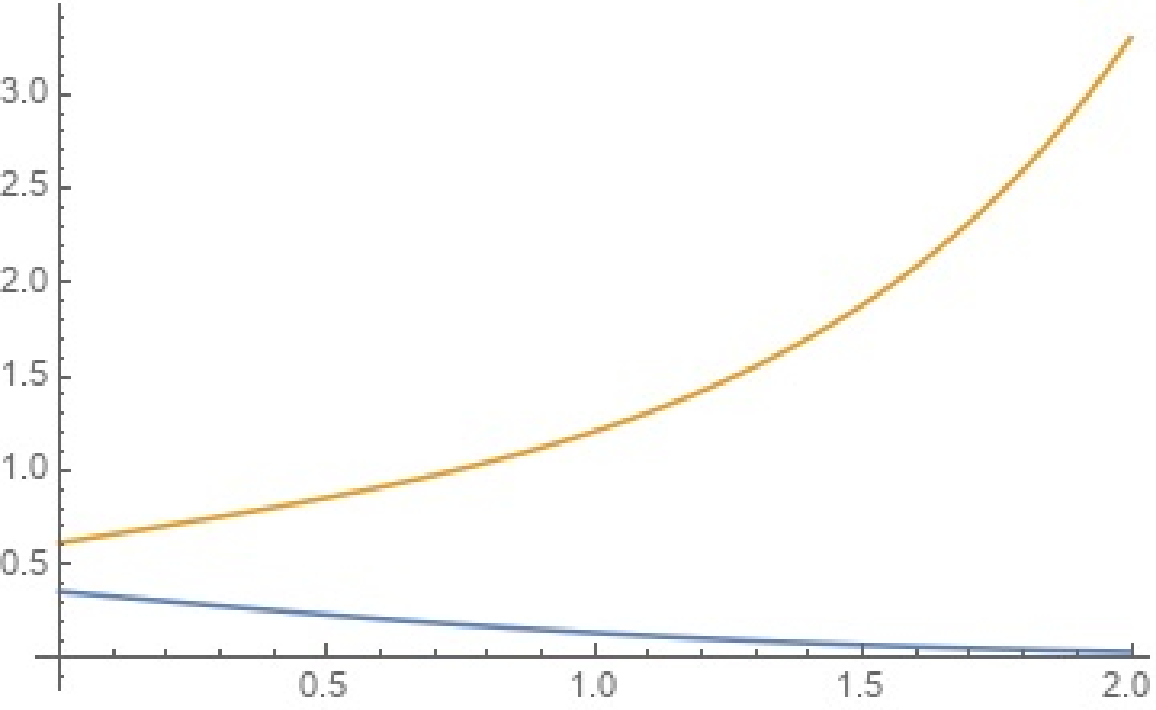}
\caption{Airy ${\rm Ai}(x)$ and ${\rm Bi}(x)$ drawn by Mathematica}
\label{fig:airy}
\end{center}
\end{figure}

We can also see that the condition 3 of the Situation \ref{situation123} 
holds by applying
the theory of singularity of ordinary differential equations
(see, e.g., the manual {\tt DEtools/formal\_sol} of Maple \cite{maple}
and its references
on the theory, which has a long history).
In fact, the general solution of the Airy differential equation is expressed as
\begin{eqnarray*}
 && C_1 t^{-1/4} \exp\left(-\frac{2}{3} t^{3/2} \right) (1+ O(t^{-3/2}))  \\
 &+&C_2 t^{-1/4} \exp\left(\frac{2}{3} t^{3/2} \right) (1+ O(t^{-3/2})) 
\end{eqnarray*}
when $t \rightarrow +\infty$
where $C_i$'s are arbitrary constants.

We note that the high precision evaluation of the Airy function is studied
by several methods (see, e.g., \cite{CM2013} and its references).
Some mathematical software systems have evaluation functions of the Airy function.
For example, {\tt N[AiryAi[10]]} gives the value of ${\rm Ai}(10)$ on Mathematica.
By utilizing these advanced evaluation methods,
we use the Airy differential equation for our test case to check the validity
of our heuristic algorithm.
\end{example}

We are going to propose some heuristic methods to avoid the instability
problem of the HGM.
Numerical schemes such as the Runge-Kutta method obtain a numerical solution
by the recurrence
\begin{equation}
 F_{k+1} = Q(k,h) F_k
\end{equation}
from $F_0$
where $Q(k,h)$\footnote{It was denoted by $Q(t_0+kh,h)$ in the previous section.
We denote $Q(t_0+kh,h)$ by $Q(k,h)$ as long as no confusion arises.}
is an $r \times r$ matrix determined by a numerical scheme
and $h$ is a small number
The vector $F_k$ is an approximate value of $F(t)$ at
$t = t_k = t_0 + hk$. 

\begin{example}\rm
The Euler method assumes 
$dF/dt(t)$ is approximated by $(F(t+h)-F(t))/h$ and the scheme of this method
is
$$ F_{k+1} = (E + h P(t_k)) F_k$$
where $E$ is the $r \times r$ identity matrix.
\end{example}

In case that the initial value vector $F_0$ contains an error,
the error may generate a blow-up solution ${\tilde F}$
under the Situation \ref{situation123}
and we cannot obtain the true solution.

Let $N$ be a suitable natural number and put
\begin{equation} \label{eq:bmatrix}
  Q=Q(N-1,h) Q(N-2,h) \cdots Q(1,h) Q(0,h)
\end{equation}
We call $Q$ the {\it matrix factorial} of $Q(k,h)$.
The matrix $Q$ approximates the fundamental solution matrix of the ODE.
We assume the eigenvalues of $Q$ are positive real and distinct to simplify
our presentation.
The following heuristic algorithm avoids to get the blow-up solution.
\begin{algorithm} \label{alg:simple}
\begin{enumerate}
\item Obtain eigenvalues $\lambda_1 > \lambda_2 > \cdots > \lambda_r > 0$ of $Q$
and the corresponding eigenvectors $v_1, \ldots, v_r$.
\item Let $\lambda_m$ be the first non-positive eigenvalue. 
\item Express the initial value vector $F_0$ containing errors in terms of $v_i$'s as
\begin{equation}  \label{eq:F0_by_vi}
  F_0 = f_1 v_1 + \cdots + f_r v_r, \quad f_i \in {\bf R}
\end{equation}
\item Choose a constant $c$ such that
$ F_0' := c(f_m v_m + \cdots + f_r v_r) $ approximates $F_0$.
\item Determine $F_N$ by $F_N = Q F_0'$ with the new initial value vector
$F_0'$.
\end{enumerate}
\end{algorithm}
We call this algorithm the {\it defusing method} for initial value problem.
This is a heuristic algorithm and
the vector $F_0'$ gives a better approximation of the initial value vector
than $F_0$ in several examples and we can avoid the blow-up of the numerical solution. 

\begin{example}\rm
We set $t_0=0$, $h=10^{-3}$, $N=10 \times 10^3$
and use the 4-th order Runge-Kutta scheme.
We have
$\lambda_1 = 9.708 \times 10^9$,
$v_1 = (-5.097,-159.919)^T$
and
$\lambda_2 = 3.247\times 10^{-7}$,
$v_2 = (-5.09798,37.164813649680576037539971418209465086)^T = (a,b)$.
Then, $m=2$.
We assume the $3$ digits accuracy of the value ${\rm Ai}(0)$ as
$0.355$ and
set $F_0' = (0.355, 0.355 b/a)$.
Then, the obtained value $F_{5000}$ at $t=5$ is
$(0.000108088745179140, -0.000246853220440734)$.
We have the following accurate value by Mathematica
\begin{verbatim}
In[1]:= N[AiryAi[5]]
Out[1]= 0.000108344 
In[2]:= N[D[AiryAi[x],{x}] /. {x->5}]
Out[2]= -0.000247414
\end{verbatim}
Note that $3$ digits accuracy has been kept for the value ${\rm Ai}(5)$.
On the other hand, we appy the 4th order Runge-Kutta method
with $h=10^{-3}$ 
for $F_0 = (0.355,-0.259)^T$,
which has the accuracy of $3$ digits.
It gives the value at $t=5$ as
$(-0.147395,-0.322215)$,
which is a completely wrong value,
and 
the value at $t=10$ as 
$(-102173,-320491)$,
which is a blow-up solution.
\end{example}

This heuristic algorithm avoids the blow-up of the numerical solution.
Moreover, when the numerical scheme gives a good approximate solution
for the exact initial value,
we can give an error estimate of the solution by our algorithm.
Let $| \cdot |$ be the Eucledian norm.
\begin{lemma}  \label{lem:error1}
Let $F(t)$ be the solution.
When $| Q \F0true - F(Nh) | < \delta$ holds,
we have
\begin{equation}
 | Q F_0' - F(Nh) | < | QF_0'| + |F(Nh)| + 2 \delta
\end{equation}
for any $F_0' \in {\bf R}^n$.
\end{lemma}

{\it Proof}\/.
It is a consequence of the triangular inequality.
In fact, we have
\begin{eqnarray*}
 & & | Q F_0' - F(Nh) | \\
 &=& | Q F_0' - Q \F0true + Q \F0true - F(Nh) | \\
 &\leq& | Q F_0' - Q \F0true| + |Q \F0true - F(Nh) | \\
 &\leq& | Q F_0' | + |Q \F0true| + \delta \\
 &\leq& | QF_0'| + |F(Nh)| + 2 \delta
\end{eqnarray*}
\qed

Under the Situation \ref{situation123}, $|F(Nh)|$ is small enough.
Then, it follows from the Lemma that $|Q F_0'|$ should be small.
In this context, we can give an error estimate of our algorithm.
However, our numerical experiments present that the algorithm shows a better behavior
than this theoretical error estimate. 
Then, we would like to classify our defusing method
as a heuristic method.

\begin{example}\rm (Solving Airy differential equation by the defusing method.)
\progweb{defusing/intro/2023-07-21-airy.rr}
Give initial values at $t=-20$ as\\
$F_0=[-0.17640612707798468959,0.89286285673647123840]$
(${\rm Ai}(-20)$ and ${\rm Ai}'(-20)$).
It keeps $3$ digits of accuracy at $t=5$ even when we start at $t=-20$,
but the accuracy is broken at $t=6$;
The defusing method gives the value $1.09 \times 10^{-5}$ at $t=6$,
but ${\rm Ai}(6)\sim 9.95 \times 10^{-6}$.
The graph in Figure \ref{fig:airy-20-30} looks nice without a blowing-up,
but we have to be careful that it is not accurate.
For example, the defusing method gives the value $5.19 \times 10^{-49}$
at $t=30$, but ${\rm Ai}(30)\sim 3.21 \times 10^{-49}$.
We need to correct errors at several values of $t$ by different methods
of evaluation for this problem.
On the other hand, it works more accurately for the function $H^k_n$.
See Example \ref{ex:Hkn-defusing}.
It is a future problem to give a practical error bound of the defusing method.
\begin{figure}[tb]
\begin{center}
\includegraphics[width=5cm]{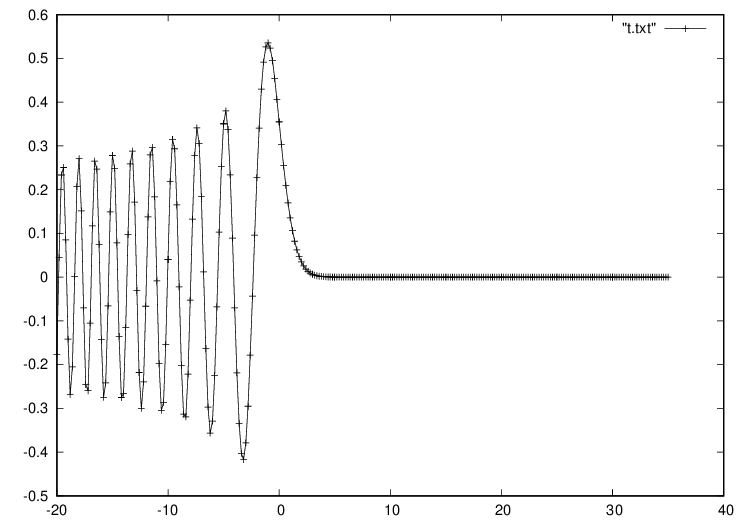}
\end{center}
\caption{Solving initial value problem of Airy differential equation by defusing, $t \in [-20,30]$}
\label{fig:airy-20-30}
\end{figure}
\end{example}

At the end of this section, we will explain how to apply defusing to the discrete QR method \cite{DRV-1994}.
The continuous orthogonormalization \cite{Davey1983} and discrete QR decomposition methods 
are based on the idea of finding a basis for the solution space 
as an orthonormal basis, 
and then determining the linear combination coefficients of these basis vectors to find the solution.

Let us explain the discrete QR decomposition method 
using the Rank 2 linear ordinary differential equation $\frac{dY}{dt} = A(t)Y$ as an example.
Assume that we numerically obtained the values of the fundamental solutions 
$Y_1, Y_2 \in \RR^2$ (column vectors) 
at $t=t_{i+1}=t_i+h$ (where $h$ is a sufficiently small number) 
from the fundamental solutions $X_1, X_2 \in \RR^2$ (column vectors) 
at $t=t_{i}$ using a difference scheme of $Y'=AY$.
In this case, $Y_1, Y_2$ are not an orthonormal basis, so we transform them into an orthonormal basis $Z_1, Z_2$ as follows.
\begin{eqnarray}
Z_1 &=& \frac{Y_1}{|Y_1|} \\
Z_2 &=& \frac{{\tilde Z}_2}{|{\tilde Z}_2|} \mbox{ where }
  {\tilde Z}_2 = Y_2 - (Z_1,Y_2)Z_1 
\end{eqnarray}
Now, we set
\begin{equation}
R(t_{i+1})=\left(\begin{array}{cc}
 |Y_1| & (Z_1,Y_2) \\
  0    & |{\tilde Z}_2| \\
\end{array}\right).
\end{equation}

Let us assume we have an initial condition $Y_0$ at $t=t_0$.
If the solution at $t=t_i$ is expressed as $c_1 X_1 + c_2 X_2$ 
using $C=(c_1,c_2)^T \in \RR^2$, 
then $c_1 Y_1 + c_2 Y_2$ is the solution at $t=t_{i+1}$.
Now, let ${\tilde C}=R(t_{i+1}) C$ and ${\tilde C}=({\tilde c}_1, {\tilde c}_2)^T$.
Then, we have
\begin{equation} \label{eq:dqr_point}
{\tilde c}_1 Z_1 + {\tilde c}_2 Z_2 = c_1 Y_1 + c_2 Y_2.
\end{equation}
Because, 
\begin{eqnarray*}
{\tilde c}_1 Z_1 + {\tilde c}_2 Z_2 &=&
 |Y_1| c_1 Z_1 + (Z_1,Y_2) c_2 Z_1 + |{\tilde Z}_2| c_2 Z_2 \\
&=& c_1 Y_1 + c_2 (Y_2-{\tilde Z}_2) + c_2 {\tilde Z}_2 \\
&=& c_1 Y_1 + c_2 Y_2.
\end{eqnarray*}
Therefore, (\ref{eq:dqr_point}) is a solution determined by the initial condition $Y_0$ at $t=t_{i+1}$.
In the Discrete QR method, the fundamental solutions $Z_1, Z_2$ are determined 
from the fundamental solutions $X_1, X_2$ at $t=t_i$, 
which form an orthonormal basis, via intermediate vectors $Y_1, Y_2$, 
to form an orthonormal basis at $t=t_{i+1}$. 
This method avoids the phenomenon of numerical explosion 
in the fundamental solutions. 
Methods for finding fundamental solutions that form an orthonormal basis are known to be stable methods for solving boundary value problems
(see, e.g., \cite{AMR}, \cite{conte-1966}, \cite{Davey1983}, \cite{DRV-1994}).

Let us apply this method to the initial value problem 
in combination with the defusing method.
Assume that the initial value at $t=t_0$ is $Y_0 \in \RR^2$.
From (\ref{eq:dqr_point}), the value of the solution at $t=t_{i+1}$ can be written as a linear transformation of $Y_0$:
\begin{equation}
(Z_1,Z_2) R(t_{i+1}) R(t_{i}) \cdots R(t_1) Y_0.
\end{equation}
Here, $(Z_1,Z_2)$ is a $2 \times 2$ matrix obtained by arranging the column vectors $Z_i$ horizontally.

Let us apply this method to the initial value problem in combination with the defusing method.
Assume the initial value at $t=t_0$ is $Y_0 \in \RR^2$.
From (\ref{eq:dqr_point}), the value of the solution at $t=t_{i+1}$ can be written as a linear transformation of $Y_0$:
\begin{equation}
(Z_1,Z_2) R(t_{i+1}) R(t_{i}) \cdots R(t_1) Y_0
\end{equation}
Here, $(Z_1,Z_2)$ is a $2 \times 2$ matrix formed by arranging the column vectors $Z_i$ horizontally.

Now, we consider a situation where $Y_0$ contains errors, 
and the subdominant solution cannot be well determined in the initial value problem.
Let us apply the defusing method
to the following matrix
\begin{equation}
(Z_1,Z_2)\prod_{k=1}^{i+1} R(t_{k}).
\end{equation}
The eigenvector component corresponding to the largest absolute value of eigenvalue is removed from $Y_0$
to find the value of the solution at $t=t_{i+1}$.
We tested this method to solve the Airy differential equation.
\progweb{orthogonal/2026-04-08-dqr-airy.rr}
In this implementation, we were able to stably obtain a subdominant solution 
(with low accuracy)
to the Airy equation.
For example, the value of the Airy function at $t=25.0108$ is
7.6885e-38 and our implementation outputs 7.65081e-38.
The value of the Airy function at $t=29.0053$ is
7.183e-47 and our implementation outputs 7.41184e-47.
However, various heuristics will be necessary depending on the problem, such as when to apply defusing and measures to prevent loss of precision.

\section{Holonomic Sparse interpolation/extrapolation methods} \label{sec:sparse-interpolation}

We found the Holonomic Sparse Interpolation/Extrapolation Method B (Holonomic SIE method B or HIE method B, abbreviated notation) 
particularly useful for analyzing large ODEs appearing in HGM using low-precision generalized boundary values obtained from Monte Carlo simulations. 
This method is a type of the least squares spectral method (LSSM) 
(see, e.g., \cite{AMR},\cite{boyd-2001}, \cite{canuto-2006}). 
The difference from conventional LSSM is that it utilizes rigorous rational number calculations or big floats as preprocessing steps for analyzing large ODEs
and ``soft constraints'' in an optimization step to avoid to overfit to noisy data.
We name Algorithm \ref{alg:sie_b} below holonomic SIE method B 
to emphasize these differences.

We solve the ODE (\ref{eq:ode})  of rank $r$ when
(approximate) values of $f(t)=Z(t)$ at 
$t=p_1, p_2, \ldots, p_{r'}$ are known.
We call the points $(p_i,q_i), q_i=f(p_i)$ {\it data points}\/.
In other words, we want to interpolate or extrapolate values of $f$
from these $r'$ values.
This is a generalization of boundary value problems
and we call this problem {\it the generalized boundary problem}\/.
The number of data points $r'$ may be more than the rank or less than the rank
of the ODE in the methods B and C in this section.

A standard numerical method to find values of $f$ on an interval
$[t_s,t_e]$ is as follows.
Devide $[t_s,t_e]$ into $N$ intervals.
Let $t_i$ be $t_s+hi$ where $h=(t_e-t_s)/N$.
We assume that the set of all $p_j$'s are a subset of the set $\{t_i\}$.
We denote by $f_i$ the value of $f(t_i)$.
We introduce the backward shift operator $\nabla f_i = f_i-f_{i-1}$.
Then
$\frac{\nabla^k f_i}{h^k}$ is approximately equal to
$f^{(k)}(t_i)$ or 
$f^{(k)}(t_{i-1})$ or
$\ldots$ or
$f^{(k)}(t_{i-k})$.
For example, when $k=1$, we have
$f(t_i)-f(t_{i-1})=f'(t_i)h + O(h^2)$ if we make the Taylor expansion
of $f(t_{i-1})=f(t_i-h)$ at $t=t_i$
and 
$f(t_i)-f(t_{i-1})=f'(t_{i-1})h + O(h^2)$ if we make the Taylor expansion
of $f(t_i)=f(t_{i-1}+h)$ at $t=t_{i-1}$.
We have
\begin{equation} \label{eq:difference_scheme}
\sum_{k=0}^r c_k(t_i) \frac{\nabla^k f_{i+s_k}}{h^k} \sim b(t_i), \quad
 0 \leq s_k \leq k.
\end{equation}
Here $s_i$ is an integer to choose an approximation of $f^{(k)}(t_i)$.
By assuming the left hand side and the right hand side are equal
and giving values $f_j$ for $t=p_j,  1 \leq j \leq r$,
we have a system of linear equations of the form
\begin{equation} \label{eq:Af=B}
 A \left(\begin{array}{c}
    f_0 \\
    f_1 \\
     \cdot \\
     \cdot \\
     \cdot \\
    f_N
  \end{array} \right)
 = B
\end{equation}
where $A$ is a $(N+1) \times (N+1)$ matrix
and $B$ is a column vector of length $N+1$.
Solving this equation, we obtain approximate values of $f(t_i)$.
We call this method {\it the sparse interpolation/extrapolation method A for HGM} in this paper.

\begin{example}\rm \label{ex:exp-airy}
We solve $Lf=0$,  \prog{intro/2021-06-11-sparse-interp.ipynb}
$ L=(\pd{t}-1)(\pd{t}^2-t)=\pd{t}^3-\pd{2}^2-t\pd{t}+t-1$
on $t \in [-9,0]$ 
with 
$f(-9)= {\rm Ai}(-9) \sim -0.0221337215473414$,
$f(-4)= {\rm Ai}(-4) \sim -0.0702655329492895$,
$f(0)= {\rm Ai}(0) \sim 0.355028053887817$
and $N=100$.
Then, we obtain approximate values of the Airy function ${\rm Ai}(t)$.
The numerical result shows that we do not have a false solution
in spite of the factor $\pd{t}-1$.
\end{example}

An alternative method to solve (\ref{eq:Af=B}) when $b=0$ 
is to construct $r$ linearly indedent solutions 
(fundamental solutions)
of (\ref{eq:difference_scheme})
and to find coefficients to express the solution of
the generalized boundary conditions as the linear combination of the 
fundamental solutions.

\bigbreak
Let us introduce 
{\it the holonomic sparse interpolation/extrapolation method B for HGM}.
This is a kind of LSSM and using it in the HGM is very useful.
We minimize 
\begin{equation}  \label{eq:square-integral}
\int_{t_s}^{t_e} |Lf(t)-b(t)|^2 d\mu(t)
\end{equation}
where $d\mu(t)$ is a measure in the $t$ space.
We approximate this integral.
A numerical integration for a function $g$ 
can be expressed as
\begin{equation}  \label{eq:numerical_integral}
I_N(g) = \sum_{j=0}^N T_j g(t_j)
\end{equation}
where $t_0=t_2 < t_1 < \cdots < t_{N-1} < t_N=t_e$
and $T_j \in {\bf R}_{\geq 0}$.
We fix a numerical integration method.
For example, the trapezoidal method can be expressed as
$h=\frac{t_e-t_s}{N}$, $t_i=t_s+hi$,
$T_j=h$ for $1\leq j < N$ and $T_0=T_N=h/2$.
We approximately expand the solution $f$ by a given basis functions
$\{e_k(t)\}$, $k=0,1, \ldots, M$
as
\begin{equation} \label{eq:basis-expansion}
f(t)=\sum_{k=0}^M f_k e_k(t), \quad f_k \in {\bf R}.
\end{equation}
We put $M+1=m$ to be compatible with the variable name in our program.
Put this expression into $Lf=b$.
We minimize the following loss function,
which is approximately equal to (\ref{eq:square-integral}),
\begin{eqnarray}
\ell(\{f_k\}) 
&=& \sum_{j=0}^N |(Lf)(t_j)-b(t_j)|^2 T_j \label{eq:lost-func} \\
&=& \sum_{j=0}^N \left|\sqrt{T_j}\sum_{k=0}^M f_k (Le_k)(t_j)-\sqrt{T_j} b(t_j) \right|^2  \nonumber
\end{eqnarray}
under the constraints\footnote{If data are given in derivative values of $f$, these constraints become relations among $q_i$'s.}  at data points
\begin{equation} \label{eq:constraint}
\sum_{k=0}^M f_k \cdot e_k(p_i)=q_i, \quad i=1,2,\ldots,r'.
\end{equation}
This is a least mean square problem under constraints
(see, e.g., \cite{nocedal}).
Since the loss function (\ref{eq:lost-func}) is defined by
the numerical integration formula (\ref{eq:numerical_integral}),
we have the following estimate.
\begin{lemma} \label{lem:error-of-method-B}
The norm $\int_{t_s}^{t_e} |Lf(t)-b(t)|^2 d\mu(t)$
is bounded by $\ell(\{f_k\})+e_N(|Lf-b|^2)$
where $e_N(|Lf-b|^2)$ is the error of the numerical integration method.
\end{lemma}
Solvers of the least mean square problem output the value
of the loss function, then we can estimate the $L^2$ norm
of $Lf-b$ by this Lemma.

We can also consider a least mean square problem with no constraints
by the loss function
\begin{equation} \label{eq:loss2}
{\tilde \ell}(\{f_k\}) = \alpha \ell(\{f_k\}) +
\beta \sum_{i=1}^r \left( \sum_{k=0}^M f_k \cdot e_k(p_i)-q_i \right)^2
+ \gamma \sum_{i=0}^N f_i^2.
\end{equation}
Here $\alpha, \beta, \gamma$ are hyperparameters for the optimization.
These hyperparameters should be adjusted according to the magnitude of the error. For example, if the error in the calculated integral value is large, 
$\beta$ should be reduced to mitigate the effect of the error.
$\gamma > 0$ is used to avoid an overfitting.
While standard Least-Squares Spectral Methods (LSSM) typically treat boundary conditions as strict constraints (setting $\beta \gg \alpha$), 
our generalized boundary values derived from Monte Carlo simulations inherently contain noise. 
Therefore, our formulation essentially adopts a machine-learning-inspired optimization framework similar to Physics-Informed Neural Networks (PINNs). 
By introducing hyperparameters and intentionally setting $\beta$ to an extremely small value (soft constraints) alongside Tikhonov regularization $\gamma$, 
we prevent high-degree polynomial basis from overfitting to the noisy data, successfully suppressing parasitic solutions.

Let us summarize the procedure of this section.
\begin{algorithm} \label{alg:sie_b}{\rmfamily
Input: data points $(p_i,q_i)$ (generalized boundary values), basis $\{e_k\}$, numerical integration scheme (\ref{eq:numerical_integral}), integration domain $[t_s,t_e]$, hyperparameters $\alpha \geq \beta, \gamma$.\\
Output: Approximation of $f$ of the form (\ref{eq:basis-expansion}).\\
Step 1: Evaluate numerically $(L\bullet e_k)(t_j)$'s, $b(t_j)$'s, $\sqrt{T_j}$,
$e_k(p_i)$'s by rational approximation of $t_j, \sqrt{T_j}, (p_i,q_i)$ and by exact rational arithmetics or big floats (Computer algebra part). \\
Step 2: Solve the least square problem (\ref{eq:loss2}) by transforming
it into the form $\|A f - B\|^2$ where $A$ is a matrix, $B$ is a vector, and 
$f=(f_0, \ldots, f_M)^T$
(the accuracy of A and B can be reduced to the accuracy required for normal numerical analysis).
}
\end{algorithm}

\begin{example}\rm \label{ex:exp-airy-B}
We solve the differential equation of Example \ref{ex:exp-airy}.
The method B does not work well on the interval $[-9,0]$
because of the oscillation of the solution,
but it works on a smaller interval $[-4,0]$.
We use Chebyshef polynomials on $[-4,0]$ upto degree $9$
as basis functions
and $h=0.01$.
Data points are 
$[[-4,-0.0702655329492895],[-3,-0.37881429],[-2,0.22740743]]$.
\prog{intro/sib-intro.py}
\prog{asir-tmp/sib-intro.rr, trybintro3()}
\end{example}

This method works well for some problems which are hard.
See the Example \ref{ex:Hkn-sib} as to an application of this method
to a 4th order ODE.
An application of this method to a problem on random matrices
is given in Example \ref{ex:ec}.
We successfully solve a rank $11$ ODE of size 25 Kbytes
appearing in this problem
by this method.
Implementations of this method for Risa/Asir, SageMath, Mathematica, Julia/OSCAR, Maple are discussed in \cite{TYZ-2026}.

\bigbreak
Let us introduce 
{\it the sparse interpolation/extrapolation method C for HGM}.
We use the basis functions $e_j(t)$ as in the method C.
We construct an $m \times m$ symmetric matrix $S$ of which $(i,j)$  element
is 
\begin{equation}
\int_{t_s}^{t_e} (Le_i(t)) (Le_j(t)) d\mu(t)
\end{equation}
Let $F$ be the column vector of length $m=M+1$:
$F=(f_0, f_1, \ldots, f_{M})^T$.
Then, we have
\begin{equation}
\int_{t_s}^{t_e} \left( L \sum_{j=0}^{M} f_j e_j(t) \right)^2 d\mu(t)
=  F^T S F.
\end{equation}
It follows from this idenity that
when $\sum_{j=0}^{M} f_j e_j(t)$ is an approximate solution of $Lf=0$,
we may expect the value of the quadratic form $F^T S F$  is small.
Let $(p_k,q_k)$'s be given data points of $(t,f(t))$.
The method C finds the vector $F$ by solving 
the quadratic programming problem
of minimizing 
\begin{eqnarray}
&&F^T S F + \sum_{k=1}^r \left(\sum_{j=0}^{M} f_j e_j(p_k)-q_k\right)^2 \nonumber \\
&=& F^T S F + \sum_{k=1}^r \left( F^T S_k^2 F-2q_kS_k F + q_k^2 \right) \\
&=& F^T S F + F^T (P_e^T P_e) F - 2 P_e^T Q + Q^T Q 
\end{eqnarray}
where $S_k={\rm diag}(e_0(p_k), \ldots, e_{M}(p_k))$
and
\begin{equation}
P_e = \left(\begin{array}{ccc}
e_0(p_1) & \cdots & e_M(p_1) \\
e_0(p_2) & \cdots & e_M(p_2) \\
\cdot & \cdots & \cdot \\
e_0(p_r) & \cdots & e_M(p_r) \\
\end{array}\right),\quad
Q=(q_1, \ldots, q_r)^T
\end{equation}
Or, we may solve the quadratic programming problem
of minimizing 
$F^T S F$ 
under the constraints
\begin{equation}
\sum_{j=0}^{M} f_j e_j(p_k)-q_k = 0, \quad k=1, \ldots, r.
\end{equation}
The method C can also solve the problem of Example \ref{ex:exp-airy-B}.
\prog{intro/sic-intro.py}
\prog{asir-tmp/sic-intro.rr, trycintro3()}


\section{The Chebyshef function method} \label{sec:chebyshef}

The approximation by Chebyshef functions is known to be very good.
We briefly overview the Chebyshef function methods for ODE's
and show that it falls in the framework of sparse interpolation/extrapolation
methods.
A comprehensive reference on the method is the book by Trefethen \cite{Trefethen}.
A nice package chebfun \cite{chebfun} for Matlab has been developed,
see also
\url{https://en.wikipedia.org/wiki/Chebfun}.
The chebfun project was initiated in 2002 by Lloyd N. Trefethen and his student Zachary Battles.

A different solver with validation and Chebyshef functions
is proposed in \cite{BBJ2018}. 
The advantage of the method is that matrices in the solver are banded
and validation is given.
We will test this method for the HGM in a next paper.

The $n$-th Chebyshef function (polynomial) is
\begin{equation}
 T_n(x) = \cos ( n \theta),  \quad x=\cos \theta 
\end{equation}
The extreme points of the curve $y=T_n(x)$ in $[-1,1]$,
which we mean points that take the values $y=1$ or $y=-1$, 
are called $n+1$ Chebyshef points (of the second kind) of $T_n$.
For example, 
$T_2(x)=2x^2-1$ and the Chebyshef points are
$\{-1,0,1\}$.

Let $f(x)$ be a function.
Fix the set of Chebyshef points for $T_n(x)$.
Let the value of $f$ at Chebyshef point $x_j$ be
$f_j$.
The (degree $n$) Chebyshev interpolant is
\begin{equation}  \label{eq:chebyshev-interpolant}
p(x) = {\sum_{j=0}^n}^\prime \frac{(-1)^j f_j}{x-x_j} /
       {\sum_{j=0}^n}^\prime \frac{(-1)^j}{x-x_j}
\end{equation}
The primes on the summation signes signify that the terms $j=0$ and $j=n$
are multiplied by $1/2$.
It takes the value $p(x_j)=f_j$ and $p(x)$ agrees with the Lagrange interpolation
polynomial for $(x_j,f_j)$'s. See, e.g., \cite[Th 5.1]{Trefethen}.

\begin{theorem} {\rm (Bernstein 1911, 1912. See, e.g., Th 8.2, Th 8.3 in \cite{Trefethen})}
If $f$ is analytic on $[-1,1]$, its Chebyshef coefficients $a_k$ 
decrease geometrically.
If $f$ is analytic and $|f|\leq M$ in the Bernstein $\rho$-ellipse about $[-1,1]$,
then $|a_k| < 2 M \rho^{-k}$.
The degree $n$ Chebyshev interpolant has accuracy
$O(M \rho^{-n})$
by the sup norm.
\end{theorem}
Here, $a_k$ is defined by
$a_k = \frac{2}{\pi}\int_{-1}^1 \frac{f(x)T_k(x)}{\sqrt{1-x^2}} dx$
and satisfies
$f = \sum_{k=0}^\infty a_k T_k(x)$
when $f$ is Lipschitz continuous.
The Bernstein $\rho$-ellipse is defined as follows.
Consider the radius $\rho$ circle in the $z$-plane. Map it by $x=(z+z^{-1})/2$
and then we obtain the $\rho$-ellipse.
Its foci is at $-1$ and $1$.
$\rho$ is the sum of semi-major and semi-minor axes.

Chapter 16 of \cite{Trefethen} shows that Chebyshev interpolants are often as good as the best approximation in practice.

Let us proceed on explaining briefly a method of solving ODE's with Chebyshef
functions and the implementation of chebfun \cite{chebfun}.
The point is an approximation of differential operators in terms of 
Chebyshef friendly matrices.

Let $X=\{X_0, X_1, \ldots, X_{n-1}\}$ be the set of the $n$ Chebyshef points (of the second kind)
for the Chebyshef function $T_{n-1}$.
The $i$-th element of $X$ $(0 \leq i < n)$ is given as
\begin{equation}
\label{eq:chebpts}
  \cos\left( \pi \frac{n-i-1}{n-1} \right). 
\end{equation}
We denote by $t$ the independent variable for ODE's instead of $x$.
We denote by $\ell_j(X;t)$ the $j$-th polynomial of the Lagrange
interpolation for $X$.
We have $\ell_j(X;X_j)=1$ and $\ell_j(X;X_i)=0$ when $i\not= j$.

When $f(t)$ is a polynomial of degree $n-1$, we have
\begin{equation}
f(t) = \sum_{j=0}^{n-1} \ell_j(X;t) f_j
\end{equation}
by putting $f_j=f(X_j)$.

Now, let $Y$ be the set of the $(n-m)$ Chebyshef points for $T_{n-m-1}$ where $m \geq 0$.
We approximate $f(t)$ by the values at $Y$,
which is called ``down-sampling'' in
\cite{DH2016}.
Since we have
\begin{equation}
f^{(m)}(t) = \sum_{j=0}^{n-1} \ell_j^{(m)}(X;t) f_j
\end{equation}
and it is a polynomial of degree $n-m-1$,
we have
\begin{eqnarray}
f^{(m)}(t) &=& \sum_{k=0}^{n-m-1} \ell_k(Y;t) 
            \left( \sum_{j=0}^{n-1} \ell_j^{(m)}(X;Y_k) f_j \right) \\
 &=& \sum_{j=0}^{n-1} \left( \sum_{k=0}^{n-m-1} \ell_k(Y;t) 
                            \ell_j^{(m)}(X;Y_k) \right) f_j
\end{eqnarray}

\begin{definition}\rm
The $(n-m) \times n$ matrix with $(i,j)$ entries
\begin{equation}
\sum_{k=0}^{n-m-1} \ell_k(Y;Y_i) \ell_j^{(s)}(X;Y_k) 
\end{equation}
is denoted by $M(n-m,n;s)$ .
\end{definition}
For $s \leq m$,
the $s$-th derivative value of $f(t)$ at $Y_i$ 
is $(\mbox{$i$-th row of $M(n-m,n;s)$}) \cdot (f_0, \ldots, f_{n-1})^T$.
Then, the matrix can be used to approximate  $\frac{d^sf}{dt^s}$.
The ODE
$$\sum_{i=0}^r c_i(t) \frac{d^i f}{dt^t} = 0$$
is approximately translated into the linear equation
\begin{equation} \label{eq:chebyshef-discretization}
 \left(\sum_{i=0}^r {\rm diag}(c_i(Y_0), c_i(Y_1), \ldots, c_i(Y_{n-r-1}))
   M(n-r,n;i) \right)
 (f_0, f_1, \ldots, f_{n-1})^T = 0
\end{equation}
in the chebfun system
where $Y$ is the $n-r$ Chebyshef points for $T_{n-r-1}$.
In the chebfun system, the matrix $M(n-m,n;s)$
is obtained by the function {\tt diffmat(n-m,n,s)}.
The solution $f(t)$ is recovered by the Chebyshef interpolant
(\ref{eq:chebyshev-interpolant}) from $f_j$'s.

The system of linear equations (\ref{eq:chebyshef-discretization})
is used in the Chebyshef function method.

\begin{example}\rm
The Airy equation 
$$
  f''-t f = 0
$$
is solved by solving a linear equation
derived by $M(n-2,n;2)$ and $M(n-2,n;0)$
with two boundary conditions.
Symbolically, we solve
\begin{equation}
  \left( M(n-2,n;2) - {\rm diag}(Y) M(n-2,n;0) \right) F = 0
\end{equation}
where $F=(f_0, f_1, \ldots, f_{n-1})^T$ with giving, e.g., values of $f_0$ and $f_{n-1}$. \\
See \url{https://www.chebfun.org/examples/ode-linear/SpectralDisc.html}.
\end{example}

\begin{figure}[tbh]
\includegraphics[width=5cm]{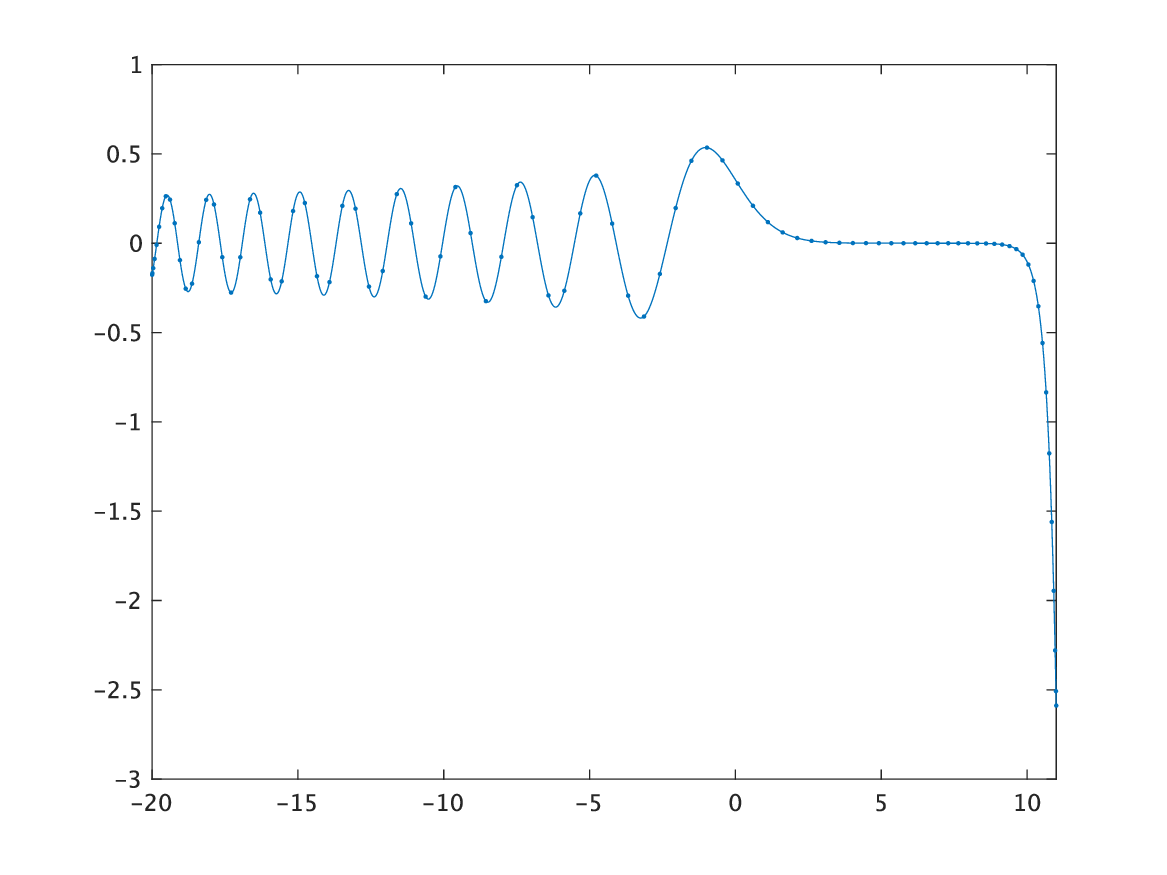}
\includegraphics[width=5cm]{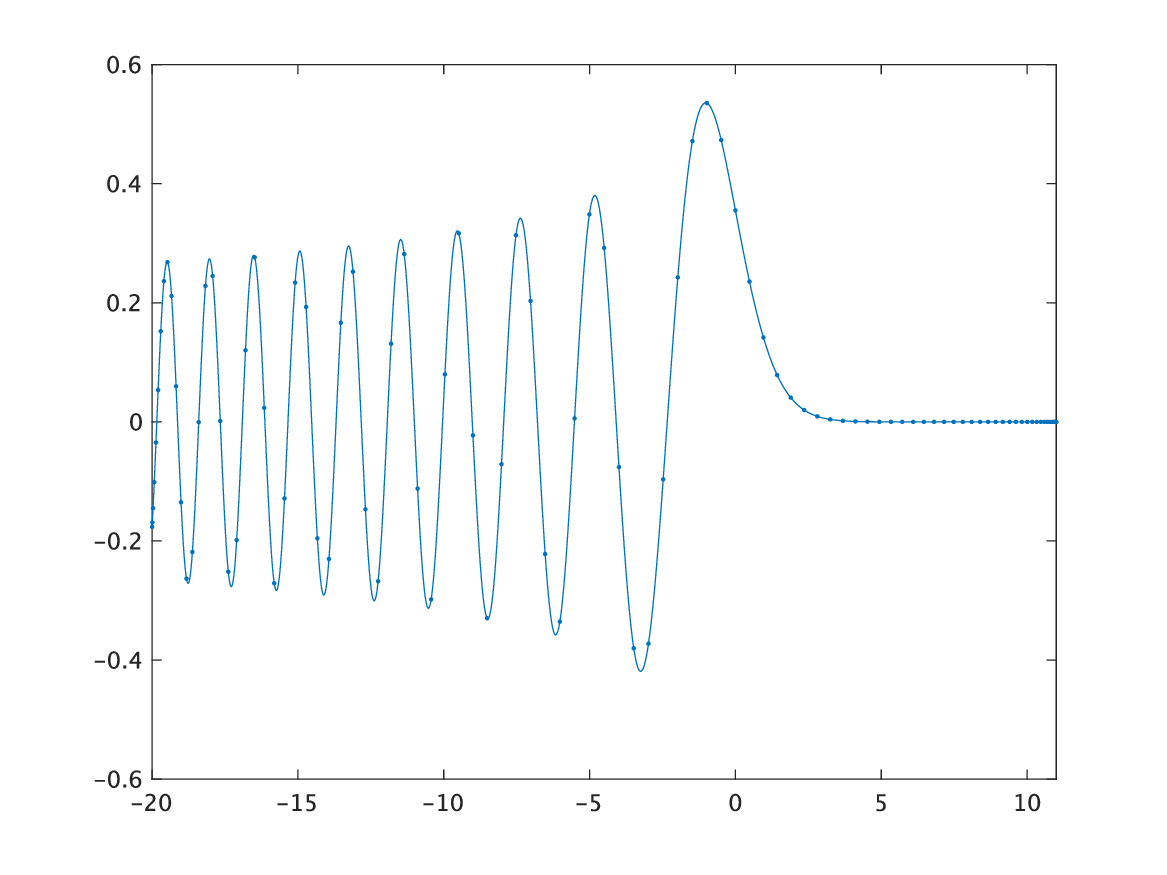}
\caption{Solving the Airy differential equation by chebfun}
\label{fig:airy-iv-bv}
\end{figure}

\begin{example}\rm
Initial value problem for Airy ${\rm Ai}(t)$.
\progweb{defusing/intro/y2023\_07\_16\_airy\_initial\_value.m}
We give the initial value at $t=-20$ as
${\rm Ai}(-20)=-0.176406127077984689590192292219$
and
${\rm Ai}'(-20)= 0.892862856736471238398409934114$.
These values are evaluated by Mathematica.
Chebfun gives reasonable values upto $t=9$,
but divergent values appear
when $t$ is larger than $9$.
See the left graph of Figure \ref{fig:airy-iv-bv}.
\end{example}

\begin{example}\rm
Boundary value problem for Airy ${\rm Ai}(t)$.
\progweb{defusing/intro/y2023\_07\_16\_airy\_boundary\_value.m}
We give the boundary value at $t=-20$ and
$t=11$ as
${\rm Ai}(-20)=-0.176406127077984689590192292219$
and
${\rm Ai}(11)=4.22627586496035959129883545080 \times 10^{-12}$.
Divergent values do not appear.
See the right graph of Figure \ref{fig:airy-iv-bv}.
\end{example}

As we have seen in two examples, the Chebyshef function
method (solving boundary value problem) 
proposed and implemented in \cite{DH2016}
works well, but solving the initial value problem is affected 
by numerical error.


Finally, we note that
the Chebyshef function method can be regarded as a special case of 
the sparse interpolation/extrapolation method B.
In fact, the numerical integration scheme of the Chebyshef quadrature
is
\begin{equation}
\int_{-1}^1 \sqrt{1-t^2} g(t) dt \sim \sum_{i=1}^{k-1} T_i g(Y_i)
\end{equation}
where $Y$ is the set of the Chebyshef points 
(\ref{eq:chebpts})
for $T_{k}$ and the weight $T_i$ is
$$
T_i = \frac{\pi}{k} \sin^2\left(\frac{i}{k}\pi\right)
$$
Put $g(t)=|Lf|^2$ and $d\mu(t) = \sqrt{1-t^2} dt$.
Since the left hand side of (\ref{eq:chebyshef-discretization}) are values
at the set of Chebyshef points $Y$,
assuming it is equal to the zero vector is equivalent to
that the integral by the Chebyshef quadrature over $Y$ is equal to zero.

\section{Tests of the methods to the function $H^k_n$}  \label{sec:Hkn}
Let $n$ and $k$ be positive integers.
We define the function $H^k_n(x,y)$ by
\begin{eqnarray}  
H^k_n(x,y) &=& \int_0^x t^k \exp(-t) {}_0F_1(;n;yt) dt \label{eq:integralHkn} \\
&=& \frac{\Gamma(n)}{\sqrt{\pi}\Gamma(n-1/2)} \int_{D(x)} t^k 
    (1-s^2)^{n-3/2} \exp(-t-2s \sqrt{yt}) dt ds  \\
 && \ \mbox{where}\     D(x)=\{(t,s) \in [0,x] \times [-1,1] \} \nonumber
\end{eqnarray}
This function appears in studies of the outage probability of MIMO WiFi systems
(see, e.g., \cite{DOTS}).
We will compare the methods in the section \ref{sec:methods}
for this function on several systems.
The function satisfies the following system of linear partial differential equations.

\begin{proposition}{\rm (\cite{DOTS})}
\begin{enumerate}
\item
The function $u=H_n^k(x,y)$ satisfies
\begin{eqnarray*}
\{\theta_{y}(\theta_{y}+n-1)+y(\theta_{x}-\theta_{y}-k-1)\}\bullet u&=&0, \\
(\theta_{x}-\theta_{y}-k-1+x)\,\theta_{x}\bullet u &=&0.
\end{eqnarray*}
where $\theta_{x} = x\frac{\pd{}}{\pd{} x}$,$\theta_{y} = y\frac{\pd{}}{\pd{} y}$.
The holonomic rank of this system is $4$.
\item The function $u$ is a solution of the following ODE with respect to $y$.
\begin{equation} \label{eq:hkn-y} 
y^2 \pd{y}^4+(-y+2n+2)y\pd{y}^3+(-y x+(-k-n-3) y+n(n+1))\pd{y}^2
+((y-n)x-n(k+2))\pd{y}+(k+1)x
\end{equation}
\end{enumerate}
\end{proposition}
\prog{Prog\_paper/19-a19-n-pf.rr}

When $y \rightarrow +\infty$, solutions of the system has
the following asymptotic behavior.
It is shown by the {\tt DEtools[formal\_sol]} function of Maple \cite{maple}. \prog{Hkn/a4-f2.ml}
\begin{eqnarray*}
h_1&=&(xy)^{-1/2(1/2+n)} \exp(-2 (xy)^{1/2}) (1 + O(1/y^{1/2})),\\
h_2&=&y^{-k-1}(1+O(1/y)), \\
h_3&=&(xy)^{-1/2(1/2+n)} \exp(2 (xy)^{1/2}) (1 + O(1/y^{1/2})), \\
h_4&=&y^{1-n+k} \exp(y) (1+O(1/y)),
\end{eqnarray*}
Which is the asymptotic behavior of the function $H^k_n(x,y)$
when $x$ is fixed?
We compare the value of $h_4$ and the value by a numerical integration
in Mathematica\footnote{The quality of {\tt HypergeometricPFQ} of mathematica
is extremely high. However, the method to evaluate hypergeometric functions
in Mathematica is still a black box. It is not easy to give a numerical
evaluator of hypergeometric functions which matches to Mathematica
in all ranges of parameters and independent variables.
}.
The integration by Mathematica is done as follows.
{\footnotesize
\begin{verbatim}
--> hh[k_,n_,x_,y_]:=NIntegrate[t^k*Exp[-t]*HypergeometricPFQ[{},{n},t*y],{t,0,x}];
--> hh[10,1,1,1000]
\end{verbatim}
}

\begin{center}
\begin{tabular}{|r|c|}
$y$ &  Ratio \\ \hline 
1000 &  7.36595030875893e-452\\
2000 &  2.64621603289928e-881\\
3000 &  2.67723893601667e-1311\\
\end{tabular}
\end{center}
where
\begin{equation} \label{eq:find_m_Hkn}
\mbox{Ratio} = (\mbox{$H^{10}_1(1/2,y)$})/(y^{1-n+k} \exp(y)).
\end{equation}
This computational experiments suggest  that
$H^k_n$ is expressed by $h_1, h_2, h_3$ without the dominant component $h_4$
as explained in (\ref{eq:find_m}).

\begin{example} \label{ex:hkn_asymp} \rm
We approximate $H^{10}_1(1,y)$ by $h_3$ of the asymptotic series terms truncated 
by $O(y^{-3})$.  \prog{Hkn/2021-06-06-asymp.rr}
Let $y_s$ be a sufficiently large number.
We determine the constant $C_s$ by setting $H^{10}_1(1,y_s)/h_3(y_s)=C_s$
and use $C_s h_3(y)$ as an approximation of $H^{10}_1(1,y)$
(a naive method).
Then we have the following data of relative errors
$r_e(y)=\frac{C_s h_3(y)-H^{10}_1(1,y)}{H^{10}_1(1,y)}$.
\begin{center}
\begin{tabular}{ccc}
$y_s$ & $r_e(y_s+10)$ & $r_e(y_s+90)$ \\ \hline
$10^2$ & $0.0763$   & $0.259$  \\ 
$10^4$ & $5.7 \times 10^{-10}$ & $5.09 \times 10^{-9}$   
\end{tabular}
\end{center}
The data tells that the naive method works
well near $y=10^4$, 
but we need to use other methods near $y=10^2$.
The sparse interpolation/extrapolation method B in the Example \ref{ex:Hkn-sib}
is a more refined variation of this method.
\end{example}

Let us apply methods presented in the Section \ref{sec:methods}
to evaluate $H^{10}_1(1,y)$.

\begin{example} \label{ex:Hkn_math_RK} \rm 
We apply the implicit Runge-Kutta method to $H^{10}_1(1,y)$ 
which is a solution of the ODE (\ref{eq:hkn-y}).
We translate the ODE into a system of ODE for $U$ by
$U=(u,u',u'',u''')^T$ and derive an ODE for
$F=U \exp(-y) y^{-(1-n+k)}$.
This is the Gauge transformation by the exponential part of $h_4$
so that the solutions keep bounded values.
Then, the column vector function $F$ satisfies $F'=PF$ where
\begin{equation}
P = \left(
\begin{array}{cccc}
\frac{  - {y}  - {k}+  {n}- 1}{ {y}}&  1& 0& 0 \\
0& \frac{  - {y}  - {k}+  {n}- 1}{ {y}}&  1& 0 \\
0& 0& \frac{  - {y}  - {k}+  {n}- 1}{ {y}}&  1 \\
\frac{  (  - {k}- 1)  {x}}{  {y}^{ 2} }& \frac{   (  - {y}+ {n})  {x}+   {n}  {k}+  2  {n}}{  {y}^{ 2} }& \frac{   {y}  {x}+   (  {k}+  {n}+ 3)  {y}  -  {n}^{ 2} - {n}}{  {y}^{ 2} }& \frac{  - {k}  - {n}- 3}{ {y}} \\
\end{array}
\right).
\label{eq:pf_for_hkn}
\end{equation}
The following data is obtained by the implicit Runge-Kutta method (Gauss method
of order 10,
\cite[IV.5]{HNW}) 
implemented in Mathematica \cite{MathematicaRK}.
The Figure \ref{fig:2021-05-16-i-RK} presents
the value of $\log H^{10}_1(1,y)$ (the almost straight graph), 
the value by the implicit Runge-Kutta method, 
and the value by {\tt NDSolve} with the {\tt ExplicitRungeKutta} option
(the graph of yellow green).  \prog{Hkn/2021-05-17-i-RK.nb}
The initial value vector is evaluated at $y=1$ by the numerical integration 
with the accuracy about 16 digits
and
the initial step size is $10^{-3}$.
The implicit Runge-Kutta method works upto about $y=25$.
This interval $[1,25]$ is larger than the valid interval $[1,4.5]$ 
of the {\tt ExplicitRungeKutta} method of {\tt NDSolve}.
\begin{figure}[tb]
\begin{center}
\includegraphics[width=8cm]{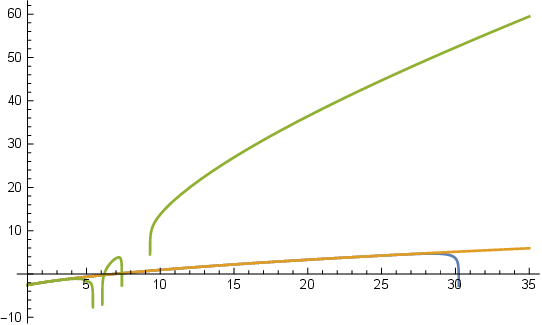}
\end{center}
\caption{The Runge-Kutta(yellow green), the implicit Runge-Kutta method(blue) starting at $y=1$, 
the exact value(orcher)} \label{fig:2021-05-16-i-RK}
\end{figure}

When we can obtain exact initial values by the numerical integration
(e.g., of (\ref{eq:integralHkn})),
it will be a good strategy to extrapolate the values in a short interval
by solving the ODE by the implicit Runge-Kutta method.
The Figure \ref{fig:2021-05-15-i-RK} describes
the value of $\log H^{10}_1(1,y)$ and the value by the implicit Runge-Kutta
method on the interval $[100,150]$.
The initial value vector is evaluated at $y=100$ by the numerical integration 
with the accuracy about 16 digits
and
the inital step size is $10^{-3}$.
It works upto about $y=130$.
\begin{figure}[tb]
\begin{minipage}{5cm}
\begin{center}
\includegraphics[width=4cm]{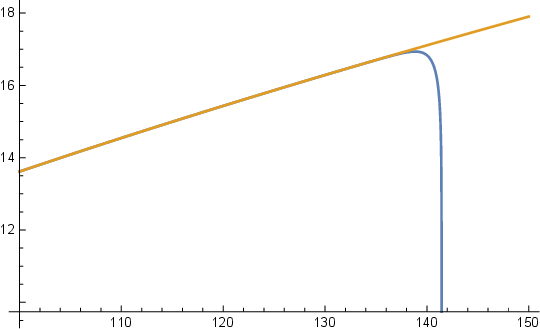}
\end{center}
\end{minipage} \quad
\begin{minipage}{5cm}
\begin{center}
\includegraphics[width=4cm]{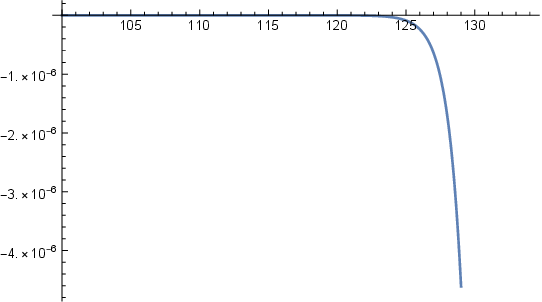}
\end{center}
\end{minipage}
\caption{The implicit Runge-Kutta method starting at $y=100$ vs the exact value.
The right graph is
the relative error $(H_r-H)/H$ where $H$ is the exact value
and $H_r$ is the value by the implicit Runge-Kutta method.
} \label{fig:2021-05-15-i-RK}
\end{figure}
\end{example}

We note that the implicit Runge-Kutta method is used 
in a large scale computation with the spectral diferred correction
(SDC), see, e.g., \cite{WSR-SDC}, \cite{pySDC}.
It is a method to accelerate the implicit Runge-Kutta method by making 
a preconditioning to iteration schemes
for solving algebraic equations of the implicit method.
We do not need this method for our relatively small rank ODE's in this paper.
However, it will be useful for solving high rank ODE's in HGM.
\prog{num-ht3/pySDC/airy\_playground.py, Airy\_implicit.py}

\begin{example} \label{ex:julia} \rm
The spectral method in the approximation theory
is a powerful method and we applied the ApproxFun package
\cite{ApproxFun} implemented in Julia for solving
the ODE (\ref{eq:hkn-y}).
The outputs are remarkable.
We solve the ODE on $[1,40]$ with the initial value at $y=1$
as $(u,u',u'',u''')(1)=(0.07810139136088563, 0.05096276584900834, 0.02050273784371611, 0.005887855153702640)$ 
where $u(y)=H^{10}_1(1,y)$.
\prog{Julia/hkn4i.jl}
Then, the output value at $y=40$ is $815.0105773595695$,
which agrees with the value of numerical integration by Mathematica 
$815.0105773587113$ upto $11$ digits.
It is also powerful to solve boundary value problems.
We give the boundary values
$(u(1),u(1.1),u(39.9),u(40))=(0.0781014, 0.0833012, 803.121, 815.011)$,
which is less accurate than the previous example.
\prog{Julia/hkn4b.jl}
Then, the output value at $y=20$ is $27.021711397385513$,
which agrees with the value of numerical integration by Mathematica
$27.021701160033859079$ upto $6$ digits.
Note that we only give the boundary values in $6$ digits accuracy.
Although it works remarkably with a proper setting, 
we have had some troubles.
For example, when we input the ODE (\ref{eq:hkn-y})$\times y^2$,
the program returns {\tt NaN}.
Another example is that it did not work for
a first order system obtained from (\ref{eq:hkn-y}).
\prog{Julia/hkn3b.jl}
\end{example}

\begin{example} \label{ex:Hkn_BDF} \rm
We can try several solvers of ODE
by {\tt solve\_ivp} on scipy/python \cite{solve_ivp}. \prog{Hkn/2021-06-02-hkn.ipynb}
We solve the initial value problem of (\ref{eq:pf_for_hkn})
on $y \in [1,30]$
with the absolute tolerance $0$ and
relative tolerance $10^{-3}$
\footnote{{\tt atol=0, rtol=1e-3}}.
Note that the default values of the absolute and relative tolerances
does not work well.
Here is a table of $y$'s for which the relative error becomes 
more than $0.3$.
Refer to \cite{solve_ivp} on details on the following methods. \prog{Hkn/2021-06-08-hkn-rerr2.ipynb}
\begin{center}
\begin{tabular}{l|c}
Method & Failure $y$ ($\mbox{relative error}>0.3$) \\ \hline
RK45  &  21.5  \\ 
Radau(Implicit Runge-Kutta) & 26.1  \\ 
BDF   &  16.8  \\ 
LSODA(Adams/BDF on ODEPACK) & 20.0  \\
\end{tabular}
\end{center}
We can also try the deSolve package in R \cite{deSolve}.
Some codes are common (e.g., FOTRAN package ODEPACK) with the python package 
and results are analogous with above. \prog{Hkn/hkn.r}
\end{example}

\begin{example} \label{ex:Hkn_python_boundary} \rm
We apply the sparse interpolation/extrapolation method A.
In other words,
we solve a generalized boundary value problem of the ODE (\ref{eq:hkn-y})
on $y \in [10^4,10^4+40]$.  \prog{Hkn/hkn-10-1-10000si-random.py}
We give generalized boundary values of $u(y)=H^k_n(1,y) \times 10^{-80}$ 
at $p_1=10^4, p_2=10^4+1333h, p_3=10^4+2000h, p_4=10^4+4000h=10^4+40$ where $h=0.01$.
We approximate derivatives as
\begin{eqnarray}
u^{(1)}(y) &=& \frac{1}{h}\left( u(y)-u(y-h) \right) \label{eq:discretize} \\
u^{(2)}(y) &=& \frac{1}{h^2}\left( u(y+h)-2u(y)+u(y-h) \right) \nonumber \\
u^{(3)}(y) &=& \frac{1}{h^3}\left( u(y+h)-3u(y) + 3u(y-h)-u(y-2h) \right) \nonumber \\
u^{(4)}(y) &=& \frac{1}{h^4}\left( u(y+2h)-4u(y+h) + 6u(y)-4u(y-h)+u(y-2h) \right) \nonumber 
\end{eqnarray}
and solve the linear equation (\ref{eq:Af=B}).
The linear equation is solved by the function {\tt linsolv}
in the scipy package \cite{scipy_linalg}.
The condition number of the matrix $A$ becomes about
$2.6 \times 10^{13}$.
We give random errors by assuming the significant digits
of the given boundary values are $3$.
The Figure \ref{fig:hkn_10_1_10000si_random} is a graph of
the solutions of $30$ tries of random errors.
It works well.
We note that the robustness with respect to random errors
depends on a choice of $p_i$'s. \prog{Hkn/hkn-10-1-10000bv-random.py}
For example, we change $p_2$ to $p_1+h$ and $p_3$ to $p_4-h$.
This stands for solving the boundary value problem
with the boundary values of $f$ and $f'$ at the boundary $p_1$ and $p_4$.
We solve the linear equation with random errors.
We can see that errors are magnified in middle
from the Figure \ref{fig:hkn_10_1_10000bv_random}.

\begin{figure}[tb]
\begin{tabular}{cc}
\begin{minipage}[t]{0.45\hsize}
\begin{center}
\includegraphics[width=5cm]{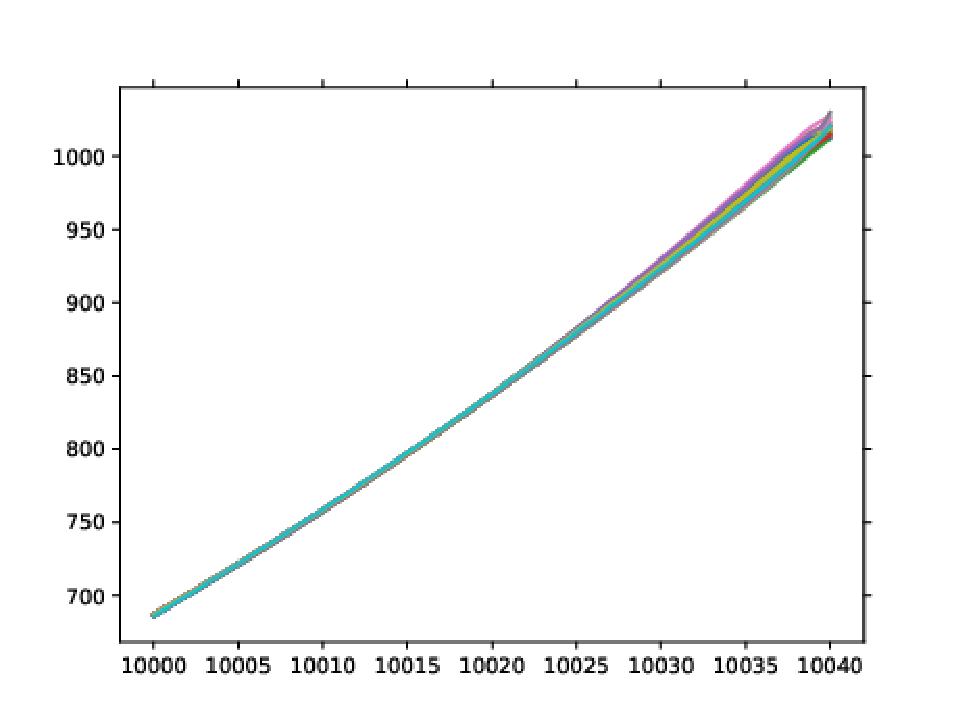}
\end{center}
\caption{Solving by the sparse interpolation A} \label{fig:hkn_10_1_10000si_random}
\end{minipage} &
\begin{minipage}[t]{0.45\hsize}
\begin{center}
\includegraphics[width=5cm]{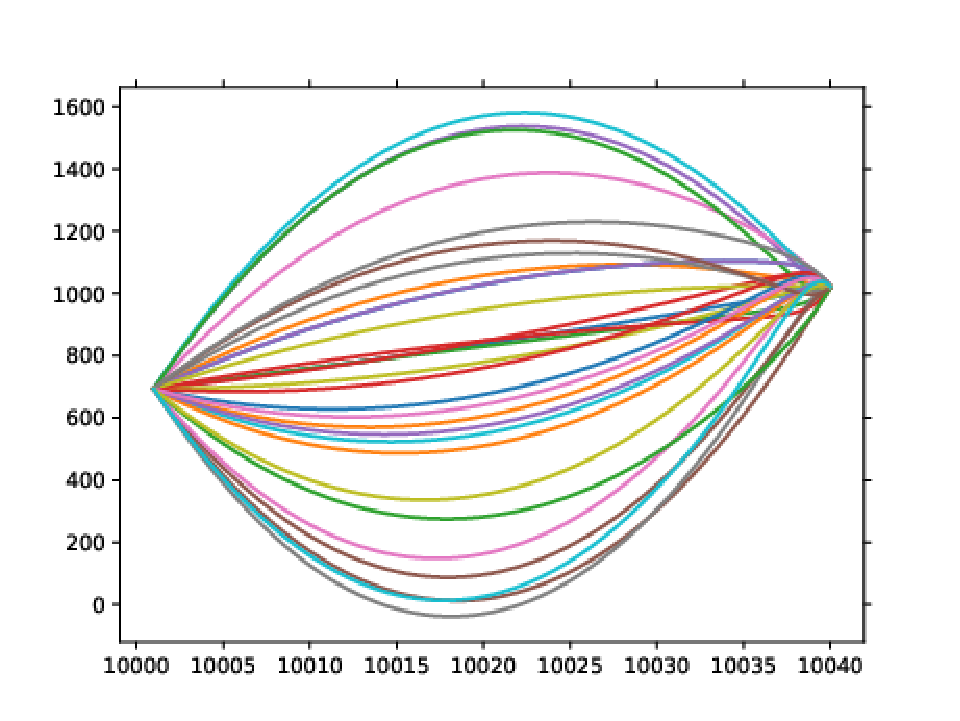}
\end{center}
\caption{Solving the boundary value problem with random errors}
\label{fig:hkn_10_1_10000bv_random}
\end{minipage}
\end{tabular}
\end{figure}

Although the boundary value problem with exact boundary values
of $H^k_n(1,y)$ and its first derivatives can be solved nicely
on $[1,40]$,
application of this method on $[1,40]$ with internal reference points
is not good. \prog{Hkn/hkn-10-1-1-40-si.py, Hkn/2021-06-11-sparse-interp.ipynb}
We take the four points as $1, 13.99675, 20.5, 40$
and $h=0.00975$ and apply the method.
It has an unacceptable relative error $33.23$ in the interval $[1,6]$,
but the relative error becomes acceptable smaller errors out of this interval.

\prog{Hkn/2021-06-09-riccati-Hkn-rtol.ipynb}
A clever method to solve boundary value problems and isolate stable subdominant solutions is to use the Riccati transformation \cite{DOR1988} 
(see also \cite[p.149]{Acton}). 
As discussed in Section \ref{sec:defusing}, this is a continuous alternative to our heuristic defusing method. We implemented this method with \texttt{solve\_ivp} on scipy/python. It works for a small interval, e.g., $y \in [1, 2]$. But, it fails on $y \in [1, 10]$. The numerical solution of the associated Riccati equation increases from $0$ to $-1 \times 10^{65}$ by the Runge-Kutta method (method=RK4 with rtol=1e-13, atol=1e-10) and the backward equation cannot be solved because of the overflow. We tried on a smaller interval $[1, 3]$ with the Adams/BDF method (method=LSODA). Evaluation of the Riccati equation did not stop in 8 minutes on the Google Colaboratory. This underscores the practical necessity of alternative robust approaches, such as the defusing method or the sparse interpolation/extrapolation methods presented in this paper, for handling the stiff instability of the HGM.
\end{example}

\begin{example} \rm \label{ex:Hkn-sib}
We apply the sparse interpolation/extrapolation method B explained 
in Section \ref{sec:sparse-interpolation}.
By utilizing the local asymptotic solution expansion at the infinity of the ODE
(\ref{eq:hkn-y}),
we use the set of the $4$ functions
\begin{equation}  \label{eq:Hkn-basis4}
  e_j(t)=t^{-3/4} \exp(2 t^{1/2}) t^{-j/2}, \quad j=0,1,2,3, \ t=y
\end{equation}
as the basis for (\ref{eq:basis-expansion}).
We will approximate the solution on $[t_s,t_e]$.
The $p_i$'s in the constraint (\ref{eq:constraint}) are
$$
p_0=t_s, p_1=t_s+5, p_2=t_s+10, \ldots, p_k=t_e-1.
$$
We do not use the constraint and 
add the difference of the approximate value and the true value 
to the the loss function
as 
$$
 \sum_{j=0}^k \left( (\mbox{approximate value at $p_j$}) -
  H^{10}_1(1,p_j) \right)^2
$$ 
We expect that this basis will give a good approximation when $t_s$ and $t_e$
are large.
In fact, when $[t_s,t_e]=[1,40]$,
the maximum of the relative error is $5.47$ by our implementation on 
{\tt least\_squares} in {\tt scipy/python}.
On the other hand, this basis gives a good approximation
on $[20,20+40]$ and $[10^4,10^4+40]$.
We give random errors to the value of $q_j=H^{10}_1(1,p_j)$ as
$q_j \times (1+O(10^{-3}))$.
It is also robust to these random errors. \prog{Hkn/hkn-10-1-10000sib-random.py Hkn/hkn-10-1-20-60sib-random.py}
\begin{center}
\begin{tabular}{l|r|r}
$[t_s,t_e]$ & $[20,60]$  & $[10^4,10^4+40]$ \\ \hline
max of relative errors with exact $q_j$ 
                        & $6.21 \times 10^{-3}$ & $2.67 \times 10^{-12}$ \\
max of relative errors with $30$ random errors
                        & $1.39 \times 10^{-2}$ & $4.07 \times 10^{-3}$ 
\end{tabular}
\end{center}
\end{example}

\begin{figure}
\begin{center}
\includegraphics[width=5cm]{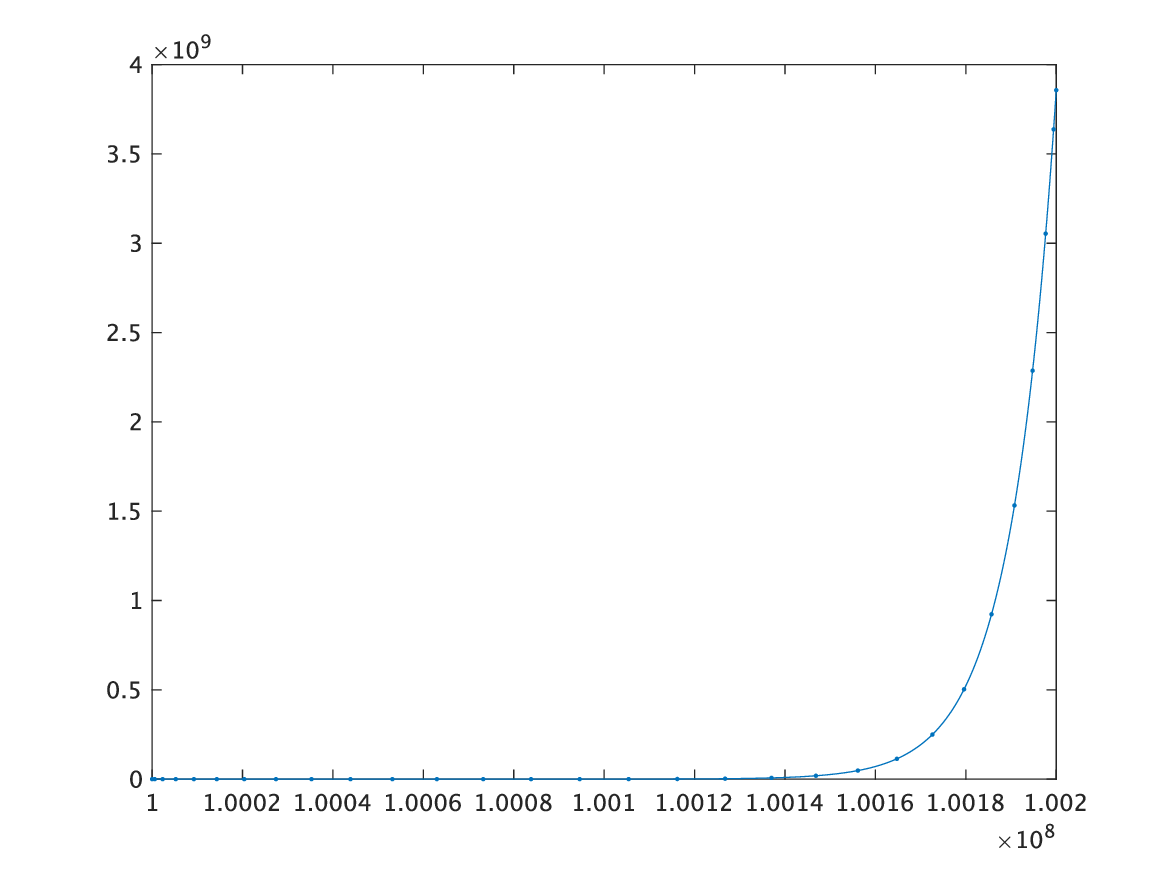}
\includegraphics[width=5cm]{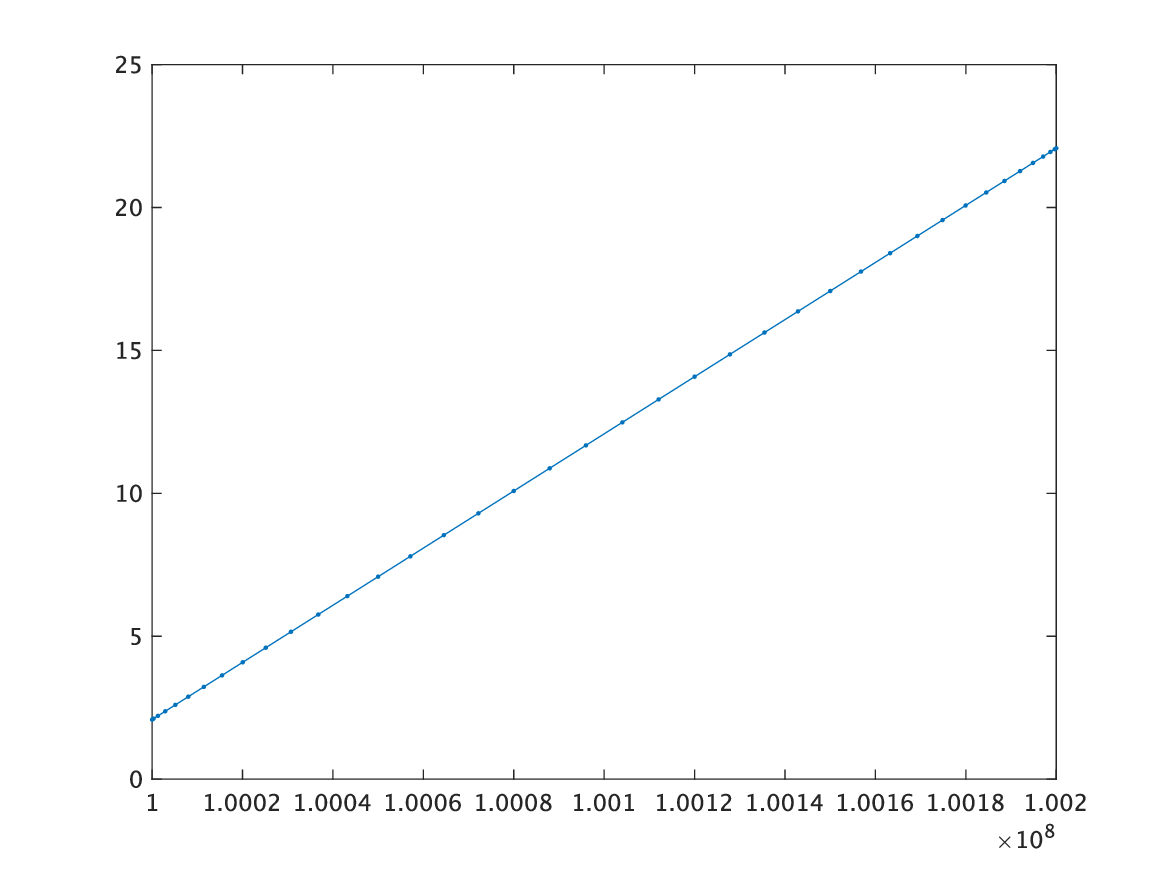}
\end{center}
\caption{Solving a boundary value problem by chebfun. Left: $H^{10}_1(1,y)$. Right: $\log H^{10}_1(1,y)$.
Values should be magnified by  $10^{8678}$.}
\label{fig:hkn-10power8}
\end{figure}
\begin{example}\rm
(Boundary value problem for $H^k_n(x,y)$ 
for $x=1$ and $y \in [10^8,10^8+2\times 10^5]$.)

We give the boundary values of $H^{10}_1(1,y)$ and 
$\frac{\partial H^{10}_1}{\partial y}(1,y)$
at $y=10^8$ and $y=10^8+2 \times 10^5$.
We apply the chebfun package for this boundary value problem.
\progweb{defusing/Hkn/y2023\_07\_25\_hkn\_valid10power8.m}
See Figure \ref{fig:hkn-10power8}.
To check the accuracy, we compare the values by the chebfun package
and by the numerical integration by Mathematica at $y=10^8+200$.
The chebfun package keeps $4$ digits accuracy at the point
and the ODE is solved in 1.66s\footnote{Apple M1, 2020, Matlab 2022b}.
On the otherhand, the numerical integration by Mathematica (2022)
\progweb{defusing/Hkn/2023-07-09-hkn-int.m}
took 23.58s\footnote{AMD EPYC 7552 48-Core Processor, 1499.534MHz}.
\end{example}

If we want to use values by the Monte-Carlo integration,
the robustness of the sparse interpolation for boundary values
with numerical errors is necessary.
\begin{figure}[tb]
\includegraphics[width=5cm]{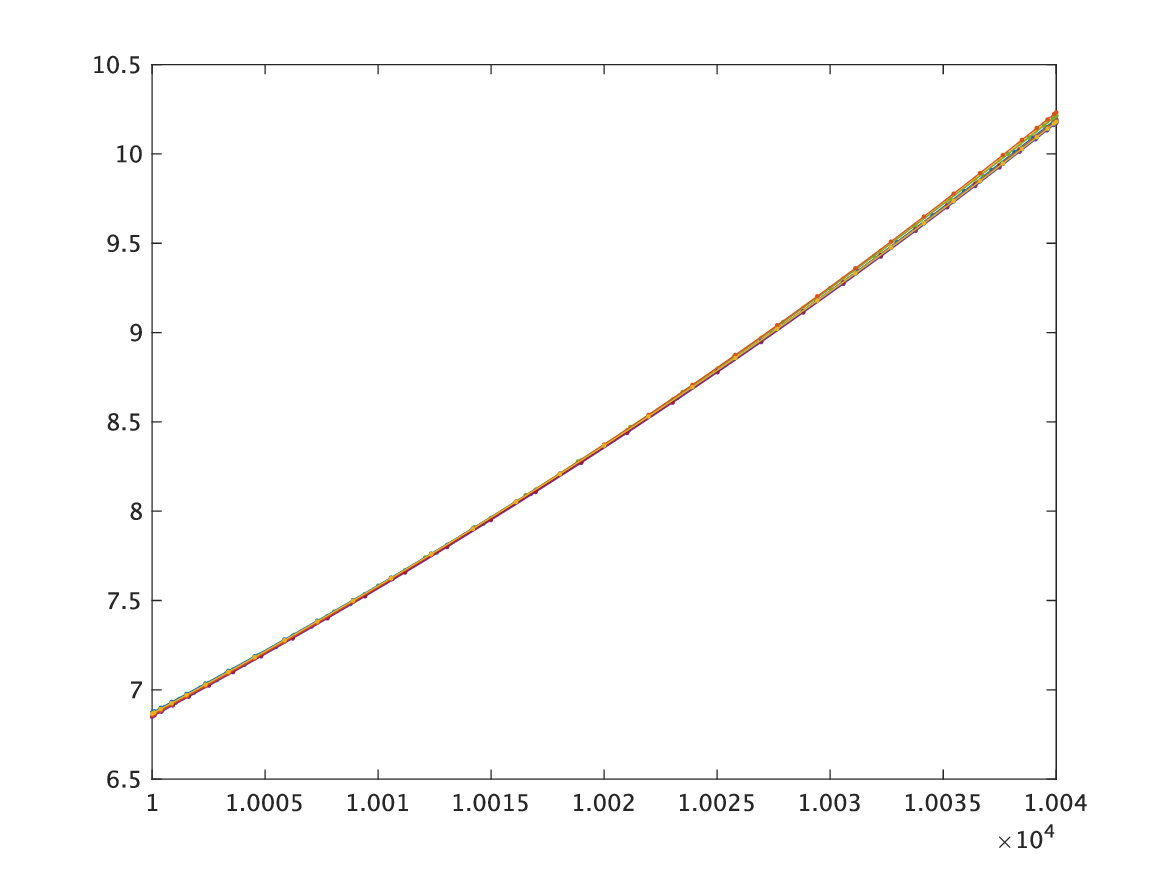}
\caption{Robustness of chebfun. $H^{10}_1(1,y)$ on $[10^4,10^4+40]$ with randomly perturbed boundary values. The values should be magnified by $10^{82}$.}
\label{fig:chebfun-robustness1}
\end{figure}
\begin{example}\rm 
(Robustness of the Chebyshev function method.)
We solve the boundary value problem for $H^{10}_1(1,y)$ on $[10^4,10^4+40]$.
We perturb the boundary values so that they keep about 3 digits accuracy.
\progweb{defusing/Hkn/y2023\_07\_27\_hkn\_10power4\_robustness.m}
Figure \ref{fig:chebfun-robustness1} is a graph of $10$ such solutions.
\end{example}

We use a gauge transformed equation to test methods other than Chebyshef 
function method. 
The chebfun system does not work well for the transformed equation.
Multiplication by the transformation seems to give some numerical errors.

In summary, we should use the implicit Runge-Kutta method 
or the sparse interpolation/extrapolation
method A when $t_s$ and $t_e$ are small and
use the sparse interpolation/extrapolation method B 
or the Chebyshef function method when they become large
for the problem $H^k_n$.
An alternative choice when $t_s$ and $t_e$ are relatively small
will be the following defusing method.

\begin{example} \rm  \label{ex:Hkn-defusing}
We implement the defusing method in {\tt tk\_ode\_by\_mpfr.rr} 
\footnote{{\tiny \url{http://www.math.kobe-u.ac.jp/OpenXM/Current/doc/asir-contrib/ja/tk_ode_by_mpfr-html/tk_ode_by_mpfr-ja.html}} (in Japanese).}
for the Risa/Asir \cite{risa-asir}.
It generates C codes utilizing the MPFR \cite{mpfr} for bigfloat and the GSL
\cite{gsl} for eigenvalues and eigenvectors.
We apply the defusing method
for initial value problem to $H^{10}_1(1,y)$ 
which is a solution of the ODE (\ref{eq:hkn-y}).
We use the step size $h=10^{-3}$ and the bigfloat of 30 digits of accuracy.
\prog{Hkn/defusing/tmp-test.c, tmp-proj.c}
\prog{Generated by asir-tmp/tk-ode-assert.rr, tk-ode-assert.hkn1(), tk-ode-assert.hkn2()}
The Figure \ref{fig:gsl-RK} shows that the adaptive Runge-Kutta method
of GSL \cite{gsl}
fails before $y$ becomes $30$. \prog{Hkn/defusing/a26-y.c}
The Figure \ref{fig:relative} presents the relative error of 
values by the defusing method and exact values.
It shows that the defusing method works even when $y=10^3$.

\begin{figure}[tb]
\begin{tabular}{cc}
\begin{minipage}[t]{0.45\hsize}
\begin{center}
\includegraphics[width=5cm]{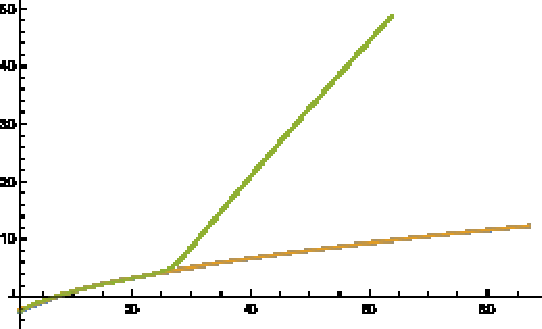}
\end{center}
\caption{$\log H^{10}_1(1,y)$. Exact value (by numerical integration) and the value by 
the defusing method agree. 
The adaptive Runge-Kutta method with the initial relative error $10^{-20}$ 
(upper curve)
does not agree with the exact value
when $y$ is larger than about $25$.}
\label{fig:gsl-RK}
\end{minipage} &
\begin{minipage}[t]{0.45\hsize}
\begin{center}
\includegraphics[width=5cm]{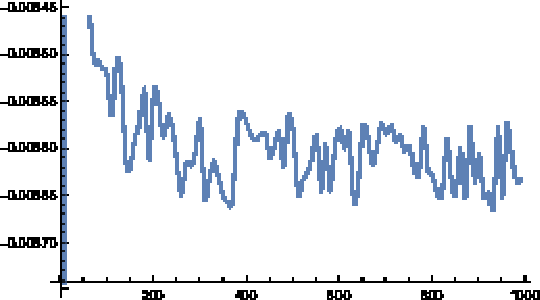}
\end{center}
\caption{The relative error of $H^{10}_1(1,y)$ of
the defusing method.
The relative error is defined as $(H_d-H)/H$
where $H_d$ is the value by the defusing method and $H$ is the exact value.
}
\label{fig:relative}
\end{minipage}
\end{tabular}
\end{figure}

\end{example}


\section{Tests of some methods to the function $E[\chi(M_t)]$} \label{sec:euler}

We consider the integral studied in \cite{TJKZ}
\begin{eqnarray}
 {\rm E}[\chi(M_t)] &=& F(s_1,s_2,m_{11},m_{21},m_{22};t) \nonumber \\
&=& \frac{1}{2}\int_t^\infty d\sigma \int_{-\infty}^\infty db \int_0^{2\pi} d\theta \int_0^{2\pi} d\phi (\sigma^2-b^2) 
\frac{s_1 s_2}{(2\pi)^2} \exp\Bigl\{-\frac{1}{2} R \Bigr\}, \label{eq:mn2}
\end{eqnarray}
where
\begin{eqnarray*}
R&= & s_{1}\left(b\sin\theta\sin\phi+\sigma\cos\theta\cos\phi-m_{11}\right)^{2}+s_{2}\left(\sigma\sin\theta\cos\phi-b\cos\theta\sin\phi-m_{21}
\right)^{2}\\ 
& & +s_{1}\left(\sigma\cos\theta\sin\phi-b\sin\theta\cos\phi\right)^{2}
 +s_{2}\left(b\cos\theta\cos\phi+\sigma\sin\theta\sin\phi-m_{22}\right)^{2}.
\end{eqnarray*}
and
$$m_{11}=1, m_{21}=2, m_{22}=3, s_1=10^3, s_2=10^2.$$
By virtue of the Euler characteristic method, 
the expectation of the Euler characteristic of a random manifold $M_t$ 
$E[\chi(M_t)]$ approximates the probability
that the maximal eigenvalue of $2 \times 2$ random matrices is less than $t$.
See \cite[Sec 3]{TJKZ} for details.

The function $E[\chi(M_t)]$  is a solution of the rank $11$ ODE
$Lf=0$
discussed in \cite[Example 5]{TJKZ}.
The operator $L$ is of the form
\begin{equation}
\left((-4.72 \times 10^{-52} t^{29}+ \cdots ) \pd{t}^{10}+ \cdots +
(-7.78 \times 10^{-22} t^{35}+ \cdots)\right) \pd{t}
\end{equation}
by multiplying a constant $10^{-31}/8.66$
\footnote{This constant is chosen so that the maximal absolute value
of the coefficients of $f_k$'s in (\ref{eq:lost-func})
of the sparse interpolation/extrapolation method B is $1$
in Example \ref{ex:ec}.}
and $\pd{t}$ from the right
to the operator given in 
\url{https://yzhang1616.github.io/ec1/ec1.html}.
\prog{or ec/tryb6/ann3.txt, variable ODE}
Note that this ODE is of a form of the singular perturbation problem 
and it seems to be hard 
to solve initial value problems like the Runge-Kutta method.
In \cite{TJKZ},
since the exact numerical integration of (\ref{eq:mn2}) is not easy,
they use values of the Monte-Carlo simulation on $t \in [3.8,3.81]$
to solve the ODE by power series.
This series approximates the expectation upto $t=3.8633$, but it does not when
$t$ is larger than it.
See \cite{TJKZ} for details.
The sparse interpolation/extrapolation method B overcomes
this difficulty.

\begin{example}\rm \label{ex:ec}
\prog{ec/tryb6/2021-07-09-try6b-tmpb.py}
\prog{It is generated by asir-tmp/sib-yi5b.rr, tryb6()}
The sparse extrapolation method B gives an approximate solution
which has larger valid domain than the solution by the power series method.
We use the basis functions
\begin{equation}
e_j(t) = (t-3.8055)^j.
\end{equation}
Let $p_i=3.8+i/1000$, $i=0,1, \ldots, 9$.
The values $q_i$ at $p_i$ are
\begin{eqnarray}
&&[0.067160, 0.065485, 0.064732, 0.063315, 0.061814, \nonumber \\
&&\ \ 0.060477, 0.059611, 0.058257, 0.057520, 0.055971]
\end{eqnarray}
by a Monte-Carlo simulation.
\prog{ec/tryb6/yiex5b.r}
$(p_i,q_i)$'s are data points for the sparse extrapolation method B.
We divide the interval $[t_s,t_e]=[3.8,4.0]$ into $200$ segments
and use the trapezoidal formula for $\ell(f)$.
The optimal solution 
\begin{eqnarray}
& & [f_0, f_1, \ldots, f_{29}] \nonumber \\
&=& [0.060145405402867516, -1.20074804549872, 9.438660716835336, \nonumber \\
& & \ -29.022737131194667, -46.606911486598264, 587.9594544508735, \ldots ]
 \nonumber 
\end{eqnarray}
can be found in about 18 seconds
\footnote{Debian 10.2 with Intel(R) Xeon(R) CPU E5-4650 2.70GHz}
by {\tt least\_squares} of the scipy package \cite{least_square} on python.
The value of the loss function $\tilde \ell$ is equal to
$3.85 \times 10^{-7}$.
The integral
\begin{equation}
\int_{t_s}^{t_e} \left( L F_{29}\right)^2 dt \mbox{ where }
F_{29}(t)=\sum_{j=0}^{29} f_j e_j(t) 
\end{equation}
is equal to $2.58 \times 10^{-4}$
and
$|F_{29}(p_i)-q_i| < 3.18 \times 10^{-4}$
for all $i$.
The left graph of Figure \ref{fig:2021_07_09_tryb6_fig} is the graph of $F_{29}(t)$
and dots are approximate values of $E(\chi(M_t))$ by the Monte-Carlo
simulation.
The data points are marked with 'x'.
The right graph of Figure \ref{fig:2021_07_09_tryb6_fig} is the relative error
of the values of $F_{29}$ and the values by the simulation.
Although the relative error is rather large,
we disagree that the method is useless.
For example, we have $F_{29}(3.935)=0.00241$,
then we can conclude the value of $t$ satisfying $f(t)=0.001$ may be
found in the domain $ t > 3.935$ assuming that the relative error
is less than $1.4$.
\begin{figure}[tb]
\begin{minipage}{6cm}
\begin{center}
\includegraphics[width=4cm]{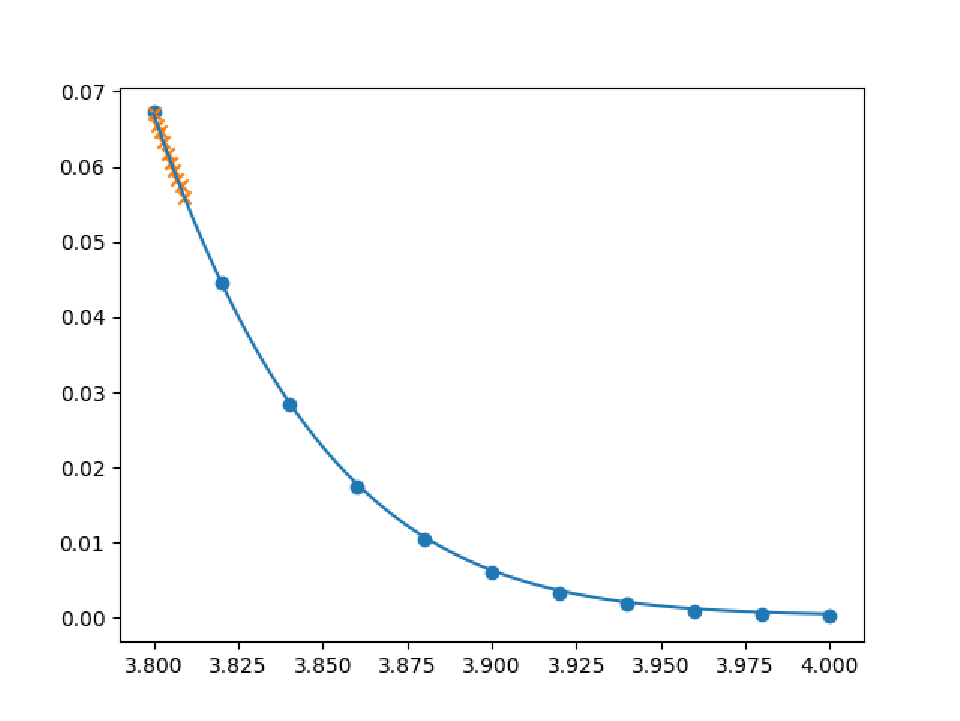}
\end{center}
\end{minipage} \quad
\begin{minipage}{6cm}
\begin{center}
\includegraphics[width=4cm]{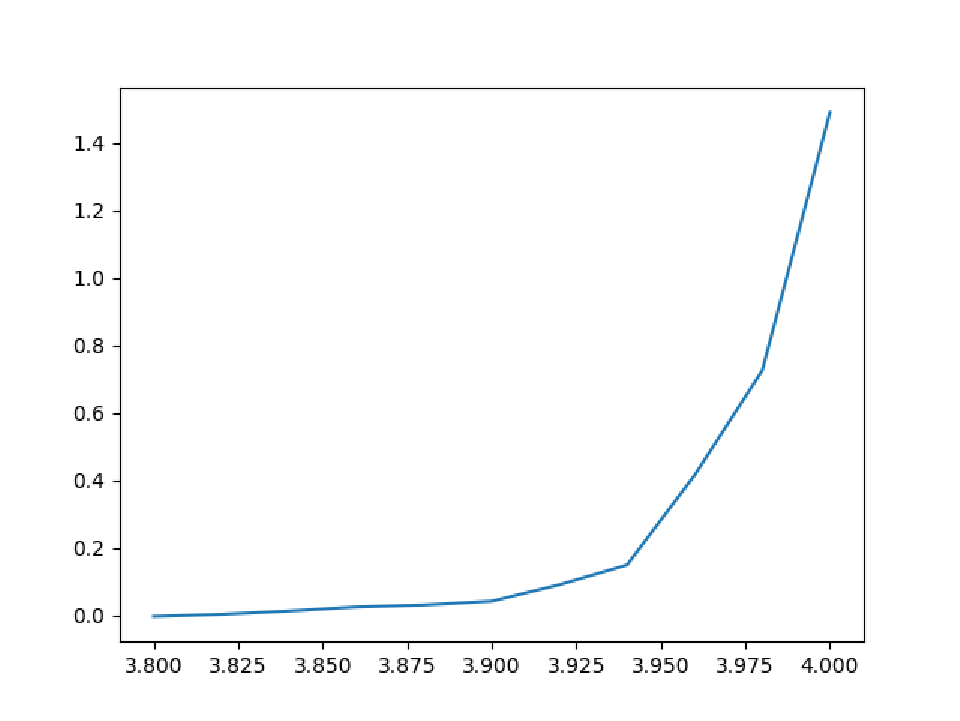}
\end{center}
\end{minipage}
\caption{
The graph of $F_{29}(t)$ and simulation values in the left and
relative errors in the right.
The data points are marked with 'x'.
} \label{fig:2021_07_09_tryb6_fig}
\end{figure}
\end{example}

\begin{example}\rm  \label{ex:ec-random}
\prog{ec/tryb5/yi5b-random.py}
We give relative errors less than $10^{-3}$ for $q_k$ of the data points $(p_k,q_k)$.
We use the same scheme as Example \ref{ex:ec}.
We execute $30$ tries and obtain the results in Figure \ref{fig:yi5b-random}.
Relative errors may be more than $40$ at far points from the data points.
We conclude that we should give accurate values of data points
to extrapolate by the sparse extrapolation method B for this problem.
\begin{figure}[tb]
\begin{minipage}{6cm}
\begin{center}
\includegraphics[width=4cm]{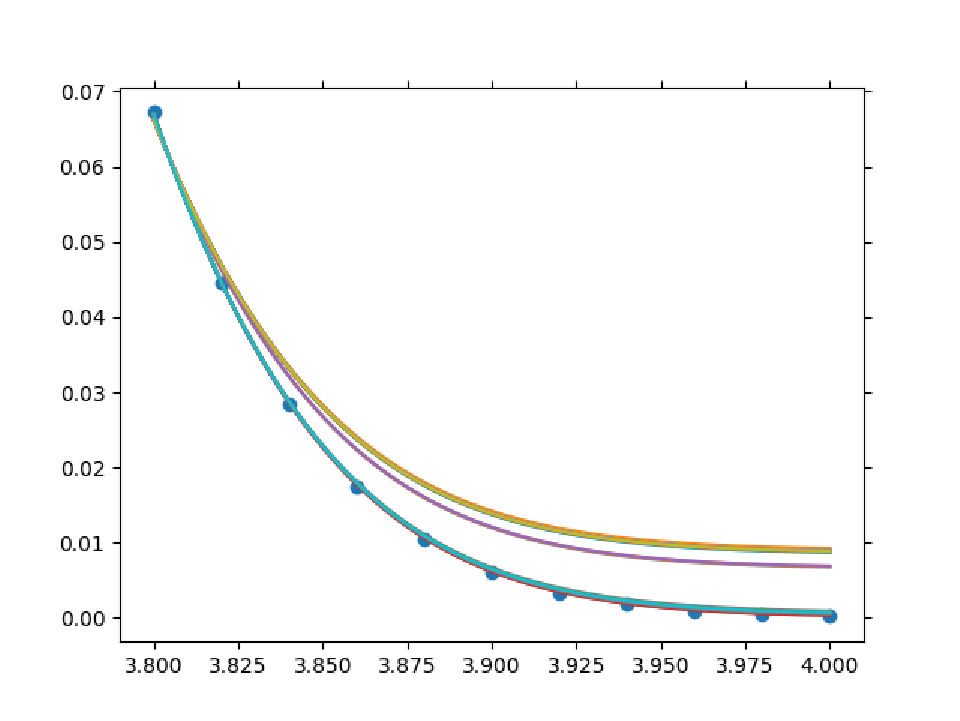}
\end{center}
\end{minipage} \quad
\begin{minipage}{6cm}
\begin{center}
\includegraphics[width=4cm]{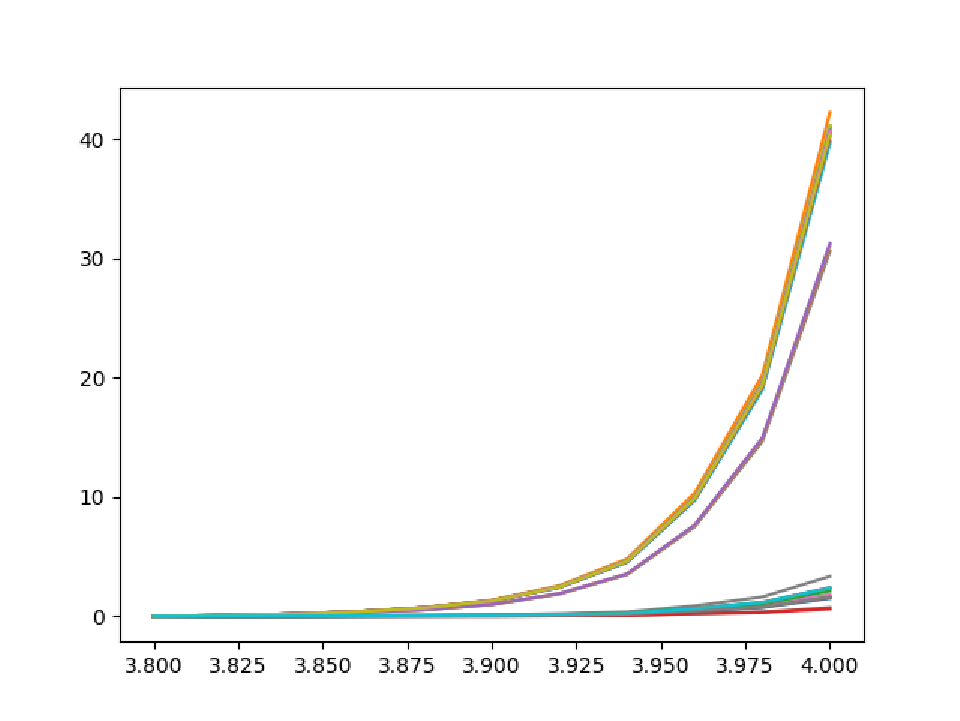}
\end{center}
\end{minipage}
\caption{
The graph of $F_{29}(t)$'s with random errors for data points and simulation values and relative errors in the right.
} \label{fig:yi5b-random}
\end{figure}
\end{example}

We expected that the sparse interpolation/extrapolation method C works
for this problem, too.
However, it does not seem to work well in real implementations on the double arithmetic.
In fact, the cvxopt package \cite{cvxopt_qp} in python returns {\tt ArithmeticError:5} in our implementation of the method C.
\prog{ec/tryb6/2021-07-14-tryc6-tmpb-modifed.py}
We expect that an implementation of the quadratic programming
with the bigfloat will solve this problem.

\bigbreak
{\it Acknowledgments}\/.
The first and the second authors are supported by the JST CREST Grant Number JP19209317.
The first author is supported by the JSPS KAKENHI Grant Number 21K03270. 
The thrid author is supported by XJTLU Research Development Funding RDF-20-01-12, NSFC Young Scientist Fund No.\ 12101506, and 
Natural Science Foundation of the Jiangsu Higher Education Institutions of China Programme-General Programme No.\ 21KJB110032.

\end{document}